\documentclass[12pt]{article}
\usepackage{amsmath,amssymb,amsfonts,fancybox,ascmac,dashbox}
\usepackage{graphicx,xcolor}
\usepackage{mathtools}
\mathtoolsset{showonlyrefs=true}
\usepackage{tikz}
\usetikzlibrary{intersections,calc}

\newcommand{\R}{\mathbb{R}}
\newcommand{\N}{\mathbb{N}}

\newcommand{\I}{\mathcal{I}}
\newcommand{\T}{\mathcal{T}}

\newcommand{\PP}{\mathcal{P}}
\newcommand{\dd}{\text{d}}
\newcommand{\lam}{\lambda}
\newcommand{\bfx}{\mathbf{x}}

\newcommand{\bfb}{\mathbf{b}}
\newcommand{\bfs}{\mathbf{s}}
\newcommand{\bfa}{\mathbf{a}}
\newcommand{\hK}{\widehat{K}}

\newcommand{\Tabc}{T_{\alpha_1\alpha_2\alpha_3}}

\newcommand{\tv}{\tilde{v}}
\newcommand{\Simp}{\mathbb S}

\newcommand{\bK}{\mathbf{K}}
\newcommand{\hT}{\widehat{T}}
\newcommand{\tT}{\widetilde{T}}
\newcommand{\bT}{\mathbf{T}}
\newcommand{\gmax}{\theta_{\max}}

\newtheorem{theorem}{Theorem}
\newtheorem{lemma}[theorem]{Lemma}
\newtheorem{corollary}[theorem]{Corollary}

\newtheorem{definition}[theorem]{Definition}
\newtheorem{assumption}[theorem]{Assumption}

\begin{document}
\title{Lectures on the Error Analysis of Interpolation \\ 
on Simplicial Triangulations\\
 without the Shape-Regularity Assumption \\[5pt]
Part 2: Lagrange Interpolation on Tetrahedrons
}
\author{Kenta Kobayashi
\footnote{Graduate School of Bussiness Administration,
         Hitotsubashi University, Kunitachi, JAPAN} \quad
Takuya Tsuchiya
\footnote{Center for Mathematical Modeling and Data Science,
Osaka University, Toyonaka, JAPAN, \newline
   \hspace*{5mm} \texttt{tsuchiya.takuya.plateau@kyudai.jp.}}}

\setlength{\baselineskip}{18pt}

\maketitle

\noindent \textbf{Abstract:}
This is the second lecture note on the error analysis of interpolation
on simplicial elements without the shape regularity assumption
\footnote{
The first one is arXiv:1908.03894 and
Memoirs of the Faculty of Science, Ehime University, 24 9--42 (2022).
}.
In this manuscript, we explain the error analysis of Lagrange
interpolation on (possibly anisotropic) tetrahedrons. This topic is hardly
explained in standard textbooks of the mathematical theory of finite
element methods. The authors hope that this manuscript will be merged 
into a new textbook in future.  Therefore, this manuscript is not 
intended to be a research paper.
Supposed readers are students and researchers who are familiar 
with the mathematical theory of the finite element methods.

\section{Lagrange interpolation on tetrahedrons}
This is the second lecture note concerning the error analysis
of interpolation on simplicial triangulations without
the shape regularity assumption.
In this note, we will explain the error analysis of Lagrange
interpolation on tetrahedrons.  To this end, we summarize the results
given in \cite{IshKobSuzTsu21, IshKobTsu20, KobayashiTsuchiya4,
 KobayashiTsuchiya5}.
Readers are referred to the first lecture note
\cite{KobayashiTsuchiya6} for the notation, lemmas, and theorems used in
this manuscript.

Throughout this paper, $T \subset \R^3$ denotes a tetrahedron
with vertices $\bfx_i$, $i=1,\cdots,4$, and all tetrahedrons are assumed to
be closed sets. Let $\lam_i$ be the barycentric coordinates of a
tetrahedron with respect to $\bfx_i$. By definition,
$0 \le \lam_i \le 1$, $\sum_{i=1}^{4} \lam_i =1$.  Let $\N_0$
be the set of nonnegative integers and
$\gamma = (a_1,\cdots,a_{4}) \in \N_0^{4}$ be a multi-index.
If $|\gamma| := \sum_{i=1}^{4}a_i = k$, then
${\gamma}/{k} := \left({a_1}/{k}, \cdots, {a_{4}}/{k}\right)$
can be regarded as a barycentric coordinate in $T$.
The set $\Sigma^k(T)$ of points on $T$ is defined by
\begin{equation*}
   \Sigma^k(T) := \left\{
   \frac{\gamma}{k} \in T \Bigm| |\gamma| = k, \;
     \gamma \in \N_0^{4} \right\}.
\end{equation*}
Let $\PP_k(T)$ be the set of polynomials defined on $T$ whose
degree is at most $k$.  For a continuous function $v \in C^0(T)$, the
Lagrange interpolation $\I_T^k v \in \PP_k(T)$ of degree $k$ is defined as
\begin{align*}
   v(\bfx) = (\I_T^k v)(\bfx), \quad \forall \bfx \in \Sigma^k(T).
\end{align*}

Let $m$, $0 \le m \le k$ be an integer, and $p$, $1 \le p \le \infty$
be a real.  For the mathematical theory of finite element methods,
estimating error $|v - \I_T^k v|_{m,p,T}$ of the Lagrange interpolation
is an important task.  For error analysis,
the following condition is usually imposed for the meshes to use in
many textbooks \cite{BrennerScott, Ciarlet, ErnGuermond}.

Suppose that $\mathcal{X}$ is a set of (possibly infinitely many)
simplicial elements (triangles or tetrahedrons).  For
$T \in \mathcal{X}$, let $h_T := \mathrm{diam}T$, and $\rho_T$ be the
diameter of its inscribed ball.
\begin{center}
\fbox{
\begin{minipage}{15truecm}
\begin{assumption}[Shape regularity]
The set $\mathcal{X}$ is called \textbf{shape regular} if there exists
a constant $\sigma > 0$ such that
\begin{align*}
  \frac{h_T}{\rho_T} \le \sigma, \qquad \forall T \in \mathcal{X}.
\end{align*}

\end{assumption}
\end{minipage}
}
\end{center}
The shape regularity assumption requires that any element
$T \in \mathcal{X}$ is not too ``flat'', or degenerate.
The maximum of the ratio $h_T/\rho_T$ in $\mathcal{X}$ is called its
\textbf{chunkiness parameter} \cite{BrennerScott}.  The shape regularity
condition is sometimes called the \textbf{inscribed ball
condition}. 

Let $\hT$ be a \textbf{reference element}.  If we consider about
tetrahedrons, the tetrahedron with vertices
$(0,0,0)^\top$, $(1,0,0)^\top$, $(0,1,0)^\top$, and $(0,0,1)^\top$ 
is typically taken as the reference element $\hT$.   Let
$\varphi(\bfx) = A\bfx + \bfb$ be an affine transformation that maps
$\hT$ to $T$, where $A$ is a $3 \times 3$ regular matrix and
$\bfb \in \R^3$. 
Error analysis is first performed on the reference element
$\hT$.  Then, the pull back $v \circ \varphi$ is used
to transfer the result obtained on $\hT$ to the ``physical
element'' $T$.  Let $\|A\|$ denote the matrix norm of $A$ associated
with the Euclidean norm of $\R^d$ ($d = 2,3$). 

Under the shape regularity assumption, we have the following theorem.
\begin{center}
\fbox{
\begin{minipage}{15truecm}
\begin{theorem}[\cite{Ciarlet}, Theorem~3.1.4]\label{Thm314}
Let $\sigma > 0$ be a constant. If $h_T/\rho_T \le \sigma$, then there
exists a constant $C = C(\hT,p,k,m)$ independent of $T$ such that,
for $v \in W^{k+1,p}(T)$,
\begin{align}
  |v - \I_T^k v|_{m,p,T} & \le C \|A\|^{k+1}\|A^{-1}\|^m
   |v|_{k+1,p,T} \notag \\
  & \le C \frac{h_T^{k+1}}{\rho_T^m} |v|_{k+1,p,T}
   \le (C\sigma^m) h_T^{k+1-m} |v|_{k+1,p,T}.
  \label{standard-est0}
\end{align}
\end{theorem}
\end{minipage}
}
\end{center}

\vspace{4mm}
If the chunkiness parameter of $\mathcal{X}$ is not small enough (say,
$\sigma > 10$), $\mathcal{X}$ is called \textbf{anisotropic}.
In numerical simulation, we sometimes need to introduce
an \textit{adaptive mesh refinement} technique.  In a process of mesh
refinements, many anisotropic elements may be generated.
With such meshes, the standard theory of finite element methods
with the shape regularity assumption cannot be applied.
The main purpose of this manuscript is to explain the error analysis
of Lagrange interpolation on tetrahedrons \textit{without} the shape
regularity assumption.

\vspace{8mm}
\begin{figure}[thb]
\begin{center}
\begin{tikzpicture}[scale=0.6]
   \coordinate (A) at (0.0,0.0);
   \coordinate (B) at (7.0,0.0);
   \coordinate (C) at (0.2,0.6);
   \draw (A) -- (B) ;
   \draw (B) -- (C) ;
   \draw (C) -- (A) ;
\end{tikzpicture}
\qquad
\begin{tikzpicture}[scale=0.6]
   \coordinate (A) at (0.0,0.0);
   \coordinate (B) at (7.0,0.0);
   \coordinate (C) at (3.5,0.6);
   \draw (A) -- (B) ;
   \draw (B) -- (C) ;
   \draw (C) -- (A) ;
\end{tikzpicture}
 \caption{Two anisotropic triangles; \textit{dagger}: the maximum angle
 is not close to $\pi$, and the circumradius is not large (left),
 and \textit{brade}: the maximum angle is
 close to $\pi$ and the circumradius is large (right).}
 \label{ani-triangle}
\end{center}
\end{figure}
\vspace{-4mm}
Let $T$ be a triangle and $R_T$ be its circumradius.
Anisotropic triangles can be categorized into only two types as depicted
in Figure~\ref{ani-triangle} (\cite{CDEFT}).
Also, as is explained in \cite{KobayashiTsuchiya6}, the ``badness'' of
an anisotropic triangle can be measured by $R_T$, and
the following theorem is known \cite{KobayashiTsuchiya6}.
\begin{center}
\fbox{
\begin{minipage}{15truecm}
\begin{theorem}[Circumradius estimates]\label{thm:circumradius-est}
Let $T$ be an arbitrary triangle.  Then, for the $k$th-order Lagrange
interpolation $\I_T^k$ on $T$, the estimation
\begin{equation} \label{eq:circumradius-est}
   |v - \I_T^k v|_{m,p,T} \le C
   \left(\frac{R_T}{h_T}\right)^m h_T^{k+1-m} |v|_{k+1,p,T}
  =  C R_T^m h_T^{k+1-2m} |v|_{k+1,p,T}
\end{equation}
holds for any $v \in W^{k+1,p}(T)$,
where the constant $C=C(k,m,p)$ is independent of the geometry of $T$.
\end{theorem}
\end{minipage}
}
\end{center}
Note that by the laws of sines, we have
\begin{align}
  \frac{R_T}{h_T} = \frac{1}{2\sin \theta_T}, \qquad
  \frac{\pi}{3} \le \theta_T < \pi
\end{align}
where $\theta_T$ is the maximum inner angle of $T$.
Hence, if there exists a constant $\gmax < \pi$ and
$\theta_T \le \gmax$, we have
\begin{align}
   |v - \I_T^k v|_{m,p,T} \le C
   \left(\frac{R_T}{h_T}\right)^m h_T^{k+1-m} |v|_{k+1,p,T}
   \le C' h_T^{k+1-m} |v|_{k+1,p,T}.
\end{align}
The condition $\theta_T \le \gmax$ is called
the \textbf{maximum angle condition} with $\gmax$
for triangles.

\vspace{4mm}
For the case of tetrahedrons, 
anisotropic tetrahedrons are usually categorized into 
nine types as depicted in Figure~\ref{ani-tetra} (\cite{CDEFT}).
Also, as we will see later, the radius of the circumsphere does
not represent the ``badness'' of an anisotropic tetrahedron.
These facts suggest that the analysis on anisotropic tetrahedrons
is much more complicated than the case of anisotropic triangles.
\vspace{3mm}
\begin{figure}[thb]
\begin{center}
\begin{tikzpicture}
\coordinate (A) at (0.0,0.0);
\coordinate (B) at (-0.2,0.1);
\coordinate (C) at (0.2,0.1);
\coordinate (D) at (0.0,3.0);
\draw (A)--(B)--(D)--(C)--(A) (A)--(D);
\draw[dashed] (B)--(C);
\end{tikzpicture}\qquad\qquad
\begin{tikzpicture}
\coordinate (A) at (0.0,0.0);
\coordinate (B) at (-0.2,1.5);
\coordinate (C) at (0.2,1.5);
\coordinate (D) at (0.0,3.0);
\draw (A)--(B)--(D)--(C)--(A) (A)--(D);
\draw[dashed] (B)--(C);
\end{tikzpicture}\qquad\qquad
\begin{tikzpicture}
\coordinate (A) at (0.0,0.0);
\coordinate (B) at (-0.2,2,0);
\coordinate (C) at (0.2,1.0);
\coordinate (D) at (0.0,3.0);
\draw (A)--(B)--(D)--(C)--(A) (B)--(C);
\draw[dashed] (A)--(D);
\end{tikzpicture}\qquad\qquad
\begin{tikzpicture}
\coordinate (A) at (0.0,0.1);
\coordinate (B) at (0.0,3.0);
\coordinate (C) at (0.3,1.5);
\coordinate (D) at (0.2,0.0);
\draw (A)--(B)--(C)--(D)--(A) (B)--(D);
\draw[dashed] (A)--(C);
\end{tikzpicture}\qquad\qquad
\begin{tikzpicture}
\coordinate (A) at (0.0,0.2);
\coordinate (B) at (0.0,3.0);
\coordinate (C) at (0.4,3.0);
\coordinate (D) at (0.4,0.0);
\draw (A)--(B)--(C)--(D)--(A) (B)--(D);
\draw[dashed] (A)--(C);
\end{tikzpicture}

\begin{tikzpicture}
\coordinate (A) at (0.0,0.0);
\coordinate (B) at (0.0,2.0);
\coordinate (C) at (1.0,1.0);
\coordinate (D) at (1.0,1.2);
\draw (A)--(B)--(D)--(A) (D)--(C)--(A);
\draw[dashed] (B)--(C);
\end{tikzpicture}\qquad
\begin{tikzpicture}
\coordinate (A) at (0.0,0.0);
\coordinate (B) at (-0.1,1.0);
\coordinate (C) at (0.0,2.0);
\coordinate (D) at (1.0,1.0);
\draw (A)--(B)--(C)--(D)--(A) (B)--(D);
\draw[dashed] (A)--(C);
\end{tikzpicture}\qquad
\begin{tikzpicture}
\coordinate (A) at (0.0,0.0);
\coordinate (B) at (1.1,1.8);
\coordinate (C) at (2.2,0.0);
\coordinate (D) at (1.1,0.7);
\draw (A)--(B)--(C)--(A) (A)--(D) (B)--(D) (C)--(D);
\end{tikzpicture}\qquad
\begin{tikzpicture}
\coordinate (A) at (0.0,0.0);
\coordinate (B) at (0.0,1.8);
\coordinate (C) at (1.8,1.8);
\coordinate (D) at (1.8,0.0);
\draw (A)--(B)--(C)--(D)--(A) (A)--(C);
\draw[dashed] (B)--(D);
\end{tikzpicture}
 \caption{Nine anisotropic tetrahedorns;
(top row from left)
\textit{spire}, \textit{spear}, \textit{spindle}, \textit{spike},
\textit{splinter}, (bottom row from left)
\textit{wedge}, \textit{spade}, \textit{cap}, \textit{sliver}.}
 \label{ani-tetra}
\end{center}
\end{figure}

K\v{r}\'{i}\v{z}ek introduced the maximum angle condition for
tetrahedrons \cite{Krizek2}.
\begin{center}
\fbox{
\begin{minipage}{15truecm}
\begin{definition}[Maximum angle condition for tetrahedrons]
Let $\gmax$, $\pi/2 \le \gmax < \pi$ be a constant.
Let $T$ be an arbitrary tetrahedron. 
If all inner angles of the faces of $T$, and all dihedral angles
between two faces of $T$ are less than or equal to $\gmax$,
$T$ is said to satisfy the \textbf{maximum angle condition} with
$\gmax$.
\end{definition}
\end{minipage}
}
\end{center}

For the error analysis of Lagrange interpolation on tetrahedrons
without the shape regularity condition, the following theorem is known
\cite{Krizek2, Duran}.
\begin{center}
\fbox{
\begin{minipage}{15truecm}
\begin{theorem} \label{thm5}
Let $\gmax$, $\pi/2 \le \gmax < \pi$ be a constant.
Suppose that a tetrahedron $T$ satisfies the maximum angle condition
with $\gmax$.  Then, there exists a constant
$C = C(\gmax, p)$ with $p > 2$ such that
\begin{equation} 
   |v - \I_T^1 v|_{1,p,T} \le C h_T |v|_{1,p,T},
   \label{esti1}
\end{equation}
where $C(\gmax,p) = \mathcal{O}((p-2)^{-1/2})$ as
$p \searrow 2$.
\end{theorem}
\end{minipage}
}
\end{center}

By this theorem, we may say that, if a tetrahedron $K$ satisfies
the maximum angle condition, the error of the linear Lagrange
interpolation is of order $\mathcal{O}(h_K)$ in $L^p$-norm
with $p > 2$. 

To extend the above estimation, a theorem similar to 
Theorem~\ref{thm:circumradius-est} was desired
\footnote{Note that Apel \cite{Apel} presents a different type of error
analysis on anisotropic meshes.}.
 For that purpose,
an immediate idea is to replace the circumradius of a triangle with the
radius of circumshpere of a tetrahedron.  However, this idea can be
immediately rejected by considering the tetrahedron $T$ with vertices
$\bfx_1 := (h,0,0)^\top$, $\bfx_2 := (-h,0,0)^\top$,
$\bfx_3 := (0,-h,h^\alpha)^\top$, $\bfx_4 := (0,h,h^\alpha)^\top$ with
$h > 0$ and $\alpha > 0$.  This tetrahedron is an example of 
\textit{sliver} (see Figure~\ref{ani-tetra}).  Setting 
$v(x,y,z) := x^2 - h^2 + h^{2-\alpha}z$, we see that
$\I_T^1 v \equiv 0$, and a simple computation yields that
$|v - \I_T^1v|_{1,\infty,T} = |v|_{1,\infty,T} \ge h^{2-\alpha}$ and
$|v|_{2,\infty,T} = 2$.  Hence, if $\alpha > 2$, an inequality such 
as the one given in Theorem~\ref{thm:circumradius-est} does not
hold for the tetrahedron, although the radius of circumshpere of the above
$T$ converges to $0$ as $h \to 0$.

To express the ``badness'' of a tetrahedron, the following definition
is given \cite{IshKobSuzTsu21, IshKobTsu20}.  Let $h_i$
$(i=1,\cdots,6)$ be the length of edges of $T$ with
$h_1 \le \cdots \le h_6 = h_T := \mathrm {diam} T$.  Then, we define
$R_T$ by
\begin{align}
   R_T := \frac{h_1h_2h_T}{|T|}h_T.
    \label{def-RT}
\end{align}
The following is the main theorem of this manuscript.
\begin{center}
 \fbox{
\begin{minipage}{15truecm}
\begin{theorem}[Main Theorem] \label{main-thm}
Let $T$ be an arbitrary tetrahedron and $R_T$ be
defined by \eqref{def-RT}.
Let $k$ and $m$ be integers with $k \ge 1$ and $0 \le m \le k$.
Let $p$ be taken as 
\begin{gather}
  \begin{cases}
      2 < p \le \infty & \text{ if } k - m = 0, \\ 
      \frac{3}{2} < p \le \infty & \text{ if } k = 1, \; m = 0,\\
      1 \le p \le \infty & \text{ if } k \ge 2 \text{ and } \; k-m \ge 1.
   \end{cases}
   \label{p-cond}
\end{gather}
For the Lagrange interpolation
$\I_T^k v$ of degree $k$ on $T$,
the following estimate holds:
\begin{align*}
  & B_p^{m,k}(T) := \sup_{u\in \T_p^{k}(T)} 
   \frac{|u|_{m,p,T}}{|u|_{k+1,p,T}} \le \, C_{k,m,p} 
  \left(\frac{R_T}{h_k}\right)^m h_T^{k+1-m}, \\
  & |v - \I_T^k v|_{m,p,T} \le  C_{k,m,p} 
    \left(\frac{R_T}{h_T}\right)^m
   h_T^{k+1-m} |v|_{k+1,p,T}, \quad
   \forall v \in W^{k+1,p}(T),
\end{align*}
where $C_{k,m,p}$ is a constant depending on $k$, $m$, and $p$.
\end{theorem}
\end {minipage}
}
\end{center}

\vspace{3mm}
\noindent
\textit{Remark.} Note that, in \eqref{p-cond} and Theorem~\ref{main-thm},
the restriction $2 < p$ for the case $k=m$ comes from the
continuity of 
the trace operator $\gamma:W^{1,p}(\bT) \ni v \mapsto v|_S \in L^1(S)$,
where $S \subset \bT$ is a non-degenerate segment (see
\cite[Section~3]{KobayashiTsuchiya4} and Lemma~\ref{lem31} in Appendix).
By the counterexamples given by Shenk \cite{Shenk} and the authors
\cite{KobayashiTsuchiya5}, we find that this restriction cannot be
improved.

For the maximum angle condition of tetrahedrons, we have the following
theorem.
\begin{center}
 \fbox{
\begin{minipage}{15truecm}
\begin{theorem}\label{thm-max-angle}
Let $T$ be an arbitrary tetrahedron and $R_T$ be
defined by \eqref{def-RT}. Then, $T$ satisfies the maximum angle
condition with $\gmax \in [\pi/2,\pi)$, if and only if
there exists a fixed constant $D = D(\gmax)$ such that
\begin{align}
    \frac{R_T}{h_T} \le D.
   \label{equiv-cond}
\end{align}
\end{theorem}
\end {minipage}
}
\end{center}
This theorem implies that, with $R_T$ given in \eqref{def-RT},
the situation for tetrahedrons is very similar to that of triangles.
We immediately obtain the following corollary.
\begin{center}
\fbox{
\begin{minipage}{15truecm}
\begin{corollary}
Let $T$ be an arbitrary tetrahedron that satisfies the maximum angle
condition with $\gmax \in [\pi/2,\pi)$.
Let $k$ and $m$ be integers with $k \ge 1$ and $0 \le m \le k$.
Let $p$ be taken as \eqref{p-cond}. For the Lagrange interpolation
$\I_T^k v$ of degree $k$ on $T$,
the following estimate holds.
\begin{align*}
   |v - \I_T^k v|_{m,p,T} \le C h_T^{k+1-m} |v|_{k+1,p,T},
   \quad \forall v \in W^{k+1,p}(T),
\end{align*}
where $C$ is a constant depending only on $k$, $m$, $p$, and
$\gmax$.
\end{corollary}
\end{minipage}
}
\end{center}

In the sequel of this lecture note, we will explain the proofs of
Theorems~\ref{main-thm}, \ref{thm-max-angle} in detail.

\section{Preliminaries}
\subsection{Notation}
A triangle with vertices $\bfx_i$ ($i=1,2,3$) is denoted
by $\triangle \bfx_1\bfx_2 \bfx_3$. The edge connecting
$\bfx_i$, $\bfx_j$ and its length are denoted by
$\overline{\bfx_i\bfx_j}$ and $|\overline{\bfx_i\bfx_j}|$,
respectively.

\subsection{The Sobolev imbedding theorem}\label{Sobolevimbedding}
Let $1 < p \le \infty$. From Sobolev's imbedding theorem and Morry's
inequality, we have the continuous imbeddings
\begin{gather*}
   W^{2,p}(T) \subset C^{1,1-3/p}(T), \quad p > 3, \\
   W^{2,3}(T) \subset W^{1,q}(T) \subset C^{0,1-3/q}(T),
      \quad \forall q > 3, \\
  W^{2,p}(T) \subset W^{1,3p/(3-p)}(T) \subset C^{0,2-3/p}(T),
      \quad \frac{3}{2} < p < 3, \\
  W^{3,3/2}(T) \subset W^{2,3}(T) \subset W^{1,q}(T)
  \subset C^{0,1-3/q}(T), \quad \forall q > 3, \\
  W^{3,p}(T) \subset W^{2,3p/(3-p)}(T) \subset W^{1,3p/(3-2p)}(T)
  \subset C^{0,3-3/p}(T), \quad 1 < p < \frac{3}{2}.
\end{gather*}
For the imbedding theorem, see \cite{AdamsFournier} and
\cite{Brezis}.  Although Morry's inequality may not be applied, the
continuous imbedding $W^{3,1}(T)$ $\subset C^{0}(T)$ still holds.
For proof of the critical imbedding, see
\cite[Theorem~4.12]{AdamsFournier} and
\cite[Lemma~4.3.4]{BrennerScott}.  
 In the following, we assume that $p$
is taken so that the imbedding $W^{k+1,p}(T) \subset C^{0}(T)$
holds, that is,
\begin{align*}
     1 \le p \le \infty, \quad  \text{ if } k+1 \ge 3
      \quad \text{ and } \quad
    \frac{3}{2} < p \le \infty, \quad \text{ if } k+1 =2.
\end{align*}

\subsection{Classification of tetrahedrons into two types}
\label{classify}
As noted in \cite{Apel, IshKobTsu20, KobayashiTsuchiya6}, 
to deal with arbitrary tetrahedrons (including anisotropic ones)
uniformly, we need to classify tetrahedrons into two types.
Let $T$ be an arbitrary tetrahedron. and  $\bfx_i$, 
$i = 1, \cdots, 4$ be its vertices.  Let $e_2$ be the shortest edge of
$T$ and $e_1$ be the longest edge connected to $e_2$.
We assume that $\bfx_1$ and $\bfx_2$ are the endpoints of $e_1$.
Let $\bfx_3$ be an endpoint of $e_2$ that is not an endpoint of $e_1$.
Then, $e_1$ and $e_2$ are edges of $\triangle\bfx_1\bfx_2\bfx_3$. 
Note that we still have two cases for assigning $\bfx_1$ and
$\bfx_2$ as the endpoints of $e_1$.

Consider the plane that is perpendicular to $e_1$ and intersects
$e_1$ at its midpoint.  Then, $\R^3$ is divided by this plane into two
half-spaces.  In this situation, we have two cases,
and tetrahedrons are classified as either Type~1 or Type~2 accordingly:
\begin{itemize}
 \item \textbf{Case~1}.
If one half-space contains three vertices and the other half-space
contains one vertex,  then $T$ is classified as Type~1.
 \item \textbf{Case~2}.
If the two half-spaces contain two vertices each,
 then $T$ is classified as Type~2.
\end{itemize}
If the plane contains a vertex, then $T$ is classified as Type~1.

We now introduce the following assignment of the vertices for each case.
\begin{itemize}
 \item If $T$ is Type~1, the endpoints of $e_2$ are $\bfx_1$ and
$\bfx_3$,  and the face $\triangle \bfx_1\bfx_3\bfx_4$
belongs to one half-space.  Let $\alpha_2 := |\overline{\bfx_1\bfx_3|}$.
 \item If $T$ is Type~2, the endpoints of $e_2$ are
$\bfx_2$ and $\bfx_3$, and $e_2$ and $\overline{\bfx_1\bfx_4}$ belongs
to the different half-spaces.  Let
$\alpha_2 := |\overline{\bfx_2\bfx_3|}$.
\end{itemize}
Define $\alpha_1 := |\overline{\bfx_1\bfx_2|}$ and 
$\alpha_3 := |\overline{\bfx_1\bfx_4|}$ for both cases.

\begin{figure}[htbp]
\centering
\begin{tikzpicture}[line width = 1pt]
   \coordinate [label=left:{$\bfx_1$}](A) at (0.0,0.0);
   \coordinate [label=right:{$\bfx_2$}](B) at (5.0,0.0);
   \coordinate [label=above:{$\bfx_4$}](D) at (1.4,3.2);
   \coordinate [label=right:{\raisebox{10pt}{$\bfx_3$}}](C) at (1.5,0.9);
   \draw (A) to node[below]{$\alpha_1$}(B);
   \draw (B) -- (D);
   \draw (A) to node[left]{$\alpha_3$} (D);
   \draw[dotted] (D) -- (C) -- (B);
   \draw[dotted] (A) to node[above]{$\alpha_2$} (C);
\end{tikzpicture}
\qquad
\begin{tikzpicture}[line width = 1pt]
   \coordinate [label=left:{$\bfx_1$}](A) at (0.0,0.0);
   \coordinate [label=right:{$\bfx_2$}](B) at (5.0,0.0);
   \coordinate [label=above:{$\bfx_4$}](D) at (1.6,3.0);
   \coordinate [label=right:{\raisebox{5pt}{$\bfx_3$}}](C) at (4.0,1.3);
   \draw (A) to node[below]{$\alpha_1$}(B);
   \draw (B) to node[right]{$\alpha_2$}(C);
   \draw (C) -- (D) -- (B);
   \draw (A) to node[left]{$\alpha_3$} (D);
   \draw[dotted] (A) -- (C);
\end{tikzpicture}
\caption{Tetrahedrons of Type~1 (left) and Type~2 (right).} \label{fig2}
\end{figure}

\subsection{Standard position of tetrahedrons}
\label{stan-posi}
For considering the geometry of tetrahedrons, it is convenient to
assign coordinates of their vertices explicitly.  Suppose that
an arbitrary tetrahedron $T$ is taken and classified as explained
in Section~\ref{classify}.   Let the parameters
$s_1$, $t_1$, $s_{21}$, $s_{22}$, $t_2$ be such that
\begin{equation}
 \begin{cases}
   s_1^2 + t_1^2  = 1, \;
  s_1 > 0, \; t_1 > 0, \quad \alpha_2 s_1 \le \frac{\alpha_1}{2}, \\
 s_{21}^2 + s_{22}^2 + t_2^2 = 1, \; t_2 > 0, \quad 
  \alpha_3 s_{21} \le \frac{\alpha_1}{2}.
 \end{cases}
 \label{eq:tetra-param}
\end{equation}

Suppose that $T$ is Type~1.  Then, using translation and
rotation, we may move $T$ as $\bfx_1 \mapsto (0,0,0)^\top$, 
$\bfx_2 \mapsto (\alpha_1,0,0)^\top$, and
$\bfx_3 \mapsto (x_3,y_3,0)^\top$ with $y_3 > 0$. 
Letting $\theta := \angle \bfx_2\bfx_1\bfx_3$ and
$s_1 := \cos\theta$, $t_1 := \sin\theta > 0$, we have
$x_3 = \alpha_2s_1$, $y_3 = \alpha_2 t_1$. Note that, by the assignment
of vertices $\bfx_i$ $(i = 1, 2, 3)$, we have $s_ 1 > 0$
(otherwise
$\overline{|\bfx_1\bfx_2|} < \overline{|\bfx_3\bfx_2|}$) and
$\alpha_2 s_1 \le \frac{\alpha_1}{2}$.
In this situation, $\bfx_4$ might be below
$xy$-plain (its $z$-coordiate is negative).  If so, we use mirror
imaging with respect to $xy$-plain to make it be above $xy$-plain (make
its $z$-coordinate positive).  Let 
$(s_{21},s_{22},t_2) := \overrightarrow{\bfx_1\bfx_4}/
|\overrightarrow{\bfx_1\bfx_4}|$.
By these procedure, we may assume without loss of generality that
$T$ of Type~1 is transformed to a tetrahedron with vertices
\begin{align}
 \bfx_1 = (0,0,0)^\top, \; \bfx_2 = (\alpha_1,0,0)^\top, \; 
    \bfx_3 = (\alpha_2 s_1, \alpha_2 t_1,0)^\top, \;
  \bfx_4 = (\alpha_3 s_{21}, \alpha_3 s_{22}, \alpha_3 t_2)^\top.
 \label{std-position1}
\end{align}
(Recall that $\alpha_2 = |\overline{\bfx_1\bfa_3}|$,
$\alpha_3 = |\overline{\bfx_1\bfa_4}|$, and
$\alpha_3s_{21} \le \alpha_1/2$ by the definition.)

If $T$ is Type~2, we may transform $T$ to a tetrahedron with
vertices
\begin{align}
 \bfx_1 = (0,0,0)^\top, \; \bfx_2 = (\alpha_1,0,0)^\top, \;
  \bfx_3 = (\alpha_1 - \alpha_2 s_1, \alpha_2 t_1,0)^\top, \;
  \bfx_4 = (\alpha_3 s_{21}, \alpha_3 s_{22}, \alpha_3 t_2)^\top,
  \label{std-position2}
\end{align}
by a similar manner. 
We refer to the coordinates in \eqref{std-position1}, 
\eqref{std-position2} as the \textbf{standard position} of $T$.
We always identify $T$ with the tetrahedron with vertices
\eqref{std-position1}, \eqref{std-position2}.   Note that we have
\begin{align}
   |T| = \frac{1}{6}\alpha_1\alpha_2\alpha_3 t_1 t_2,
   \label{T-volume}
\end{align}
where $|T|$ is the volume of $T$.

\subsection{Reference tetrahedrons}
Because we have two types of tetrahedrons, it is convenient to
introduce two reference tetrahedrons to deal\ with them uniformly.
Let $\hT$ and $\tT$ be tetrahedrons that have the 
following vertices (see Figure~\ref{ref_tetra}):
\begin{align*}
  \hT \text{ has the vertices } \; (0,0,0)^\top, \; (1,0,0)^\top, \;
  (0,1,0)^\top, \; (0,0,1)^\top, \\
  \tT \text{ has the vertices } \; (0,0,0)^\top, \; (1,0,0)^\top, \;
  (1,1,0)^\top, \; (0,0,1)^\top.
\end{align*}
\begin{figure}[tbhp]
\centering
\begin{tikzpicture}[line width = 1pt]
   \coordinate (A) at (0.0,0.0);
   \coordinate [label=above:{$z$}](B) at (0.0,3.0);
   \coordinate [label=right:{$y$}](C) at (3.0, 0.0);
   \coordinate [label=left:{$x$}] (D) at (-1.5,-0.8);
   \coordinate [label=right:{{\small$1$}}](E) at ($(A)!0.95!(B)$);
   \coordinate [label=above:{{\small$1$}}](F) at ($(A)!0.95!(C)$);
   \coordinate [label=below:{{\small$1$}}](G) at ($(A)!0.9!(D)$);
   \draw [line width=0.3pt] (E) -- (B);
   \draw [line width=0.3pt] (F) -- (C);
   \draw [line width=0.3pt] (G) -- (D);
   \draw [dotted] (E) -- (A) -- (F);
   \draw [dotted] (A) -- (G);
   \draw (E) -- (F) -- (G) -- (E);
\end{tikzpicture}
\qquad
\begin{tikzpicture}[line width = 1pt]
   \coordinate (A) at (0.0,0.0);
   \coordinate [label=above:{$z$}](B) at (0.0,3.0);
   \coordinate [label=right:{$y$}](C) at (3.0, 0.0);
   \coordinate [label=left:{$x$}] (D) at (-1.5,-0.8);
   \coordinate [label=below:{{\small$(1,1)$}}] (X) at (1.5, -0.72);
   \coordinate (Y) at ($(A)!0.41!(C)$);
   \coordinate [label=right:{{\small$1$}}](E) at ($(A)!0.95!(B)$);
   \coordinate [label=above:{{\small$1$}}](F) at ($(A)!0.95!(C)$);
   \coordinate [label=below:{{\small$1$}}](G) at ($(A)!0.9!(D)$);
   \draw [line width=0.3pt] (E) -- (B);
   \draw [line width=0.3pt] (Y) -- (C);
   \draw [line width=0.3pt] (G) -- (D);
   \draw [dotted] (E) -- (A) -- (G);
   \draw [dotted] (A) -- (X);
   \begin{scope}[line width = 0.4pt]
     \draw [dotted] (Y) -- (A);
     \draw [dotted] (X) -- (F);
   \end{scope}
   \draw (E) -- (X) -- (G) -- (E);
\end{tikzpicture}
\caption{The reference tetrahedrons $\hT$ (left) and $\tT$ (right).}
\label{ref_tetra}
\end{figure}

These tetrahedrons are called the \textbf{reference tetrahedrons}.
In the following, $\hT$ corresponds to tetrahedrons of Type~1 and
$\tT$ corresponds tetrahedrons of Type~2.  We denote the reference
tetrahedrons by $\bT$, that is, $\bT$ is either of $\{\hT,\tT\}$.

\subsection{Linear transformations}
For an arbitrary tetrahedron $T$ written as \eqref{std-position1} or
\eqref{std-position2} with parameters \eqref{eq:tetra-param}, we
consider an affine transformation from the reference tetrahedrons.
Define the matrices $\widehat{A}$, $\widetilde{A}$,
$D_{\alpha_1\alpha_2\alpha_3} \in GL(3,\R)$ by
\begin{align}
   \widehat{A} := \begin{pmatrix}
        1  & s_{1} & s_{21} \\
        0  & t_{1} & s_{22} \\
        0  & 0     & t_2    
       \end{pmatrix}, \quad
   \widetilde{A} := \begin{pmatrix}
        1  & - s_{1} & s_{21} \\
        0  & t_{1} & s_{22} \\
        0  & 0     & t_2    
       \end{pmatrix}, \quad
  D_{\alpha_1\alpha_2\alpha_3} := \begin{pmatrix}
        \alpha_1  & 0 & 0 \\
        0  & \alpha_2 & 0 \\
        0  & 0     & \alpha_3
       \end{pmatrix}.
    \label{diag-mat}
\end{align}
We immediately confirm that the following lemma holds.

\vspace{1mm}
\begin{center}
 \fbox{
\begin{minipage}{15truecm}
\begin{lemma}[\cite{KobayashiTsuchiya5}] 
Let $T$ be an arbitrary tetrahedron in the standard position
\eqref{std-position1} or \eqref{std-position2} with parameters
\eqref{eq:tetra-param}.  Then, $T$ is transformed from
the reference tetrahedron $\bT$ by
$T = \widehat{A}D_{\alpha_1\alpha_2\alpha_3}(\hT)$ for Type~1,
or $T = \widetilde{A}D_{\alpha_1\alpha_2\alpha_3}(\tT)$ for Type~2. 
\end{lemma}
\end {minipage}
}
\end{center}
The linear transformation defined by $D_{\alpha_1\alpha_2\alpha_3}$
is called the \textbf{squeezing transformation} 
\cite{KobayashiTsuchiya6}, and we will show that the squeezing
transformation does not reduce approximation property of Lagrange
interpolation at all (see Theorem~\ref{squeezingtheorem}).

\vspace{1mm}
 Note that $\widehat{A}$ and
$\widetilde{A}$ are decomposed as $\widehat{A} = X\widehat{Y}$ and
$\widetilde{A} = X\widetilde{Y}$ with
\begin{align*}
  X := \begin{pmatrix}
        1 & 0 & s_{21} \\
	0 & 1 & s_{22} \\
        0 & 0 & t_2
       \end{pmatrix}, \qquad
  \widehat{Y} := \begin{pmatrix}
        1 & s_1 & 0 \\
	0 & t_1 & 0 \\
        0 & 0 & 1
       \end{pmatrix}, \qquad
  \widetilde{Y} := \begin{pmatrix}
        1 & -s_1 & 0 \\
	0 & t_1 & 0 \\
        0 & 0 & 1
       \end{pmatrix},
\end{align*}
respectively.  We consider the singular values of $\widehat{A}$,
$\widetilde{A}$, $X$, $\widehat{Y}$, and $\widetilde{Y}$.
A straightforward computation yields 
\begin{gather*}
  \det\left(X^\top X - \mu I\right) 
    = (1 - \mu)\left(\mu^2 - 2\mu + t_2^2\right), \\
  \det\left(\widehat{Y}^\top \widehat{Y} - \mu I\right) 
  = \det\left(\widetilde{Y}^\top \widetilde{Y} - \mu I\right) 
    = (1 - \mu)\left(\mu^2 - 2\mu + t_1^2\right).
\end{gather*}
Thus, we find that, setting 
$\bfs_1 := |s_1|$ and $\bfs_2 := (s_{21}^2 + s_{22}^2)^{1/2}$,
\begin{gather}
  \|X\| = (1 + \bfs_2)^{1/2}, \quad
  \|X^{-1}\| = (1 - \bfs_2)^{-1/2}, \notag \\
  \|Y\| = (1 + \bfs_1)^{1/2}, \quad
  \|Y^{-1}\| = (1 - \bfs_1)^{-1/2}, \qquad 
   Y = \widehat{Y} \text{ or } Y = \widetilde{Y}, \notag \\
  \|A\| \le \prod_{i=1}^2(1 + \bfs_i)^{1/2}, \quad
  \|A^{-1}\| \le \prod_{i=1}^2(1 - \bfs_i)^{-1/2}, \qquad
  A = \widehat{A} \text{ or } A = \widetilde{A}.
  \label{mat-norm}
\end{gather}
Note that
\begin{align}
   \bfs_i^2 + t_i^2 = 1, \; i = 1, 2 \quad \text{ and } \quad
 \|A^{-1}\| \le \prod_{i=1}^2(1 - \bfs_i)^{-1/2} =
    \prod_{i=1}^2 \frac{(1 + \bfs_i)^{1/2}}{t_i}.
   \label{mat-norm2}
\end{align}

\subsection{Another geometric quantities of tetrahedrons}
In \eqref{def-RT}, a quantity $R_T$ is defined for a tetrahedron $T$.
Here, we define another quantity $H_T$ \cite{IshKobTsu20}, which
represent the geometry of $T$, by
\begin{align}
   H_T := \frac{\alpha_1 \alpha_2\alpha_3}{|T|}h_T
        = \frac{6h_T}{t_1t_2},
    \label{def-RT2}
\end{align}
where the last equation is from \eqref{T-volume}.
Then, the following lemma holds \cite[Lemma~3]{IshKobTsu20}.
\begin{center}
\fbox{
\begin{minipage}{15truecm}
\begin{lemma}
The two quantities $R_T$ and $H_T$ are equivalent.  That is, for
an arbitrary tetrahedron $T$, we have
\vspace{-3mm}
\begin{align}
   \frac{1}{2}H_T \le R_T \le 2 H_T.
   \label{equiv}
\end{align}
\end{lemma}
\end{minipage}
}
\end{center}
\noindent 
\textit{Proof.} Suppose that we have a triangle with the edge lengths
$h_1 \le h_2 \le h_3$.  Then, $\frac{1}{2}h_3 < h_2 \le h_3$.
Let $T$ be an arbitrary tetrahedron $T$ in the standard position.

\noindent
\textbf{Case~1.} Suppose that $T$ is of Type~1.  Set
$\beta := |\overline{\bfx_2\bfx_3}|$,
$\gamma := |\overline{\bfx_3\bfx_4}|$, and
$\delta := |\overline{\bfx_2\bfx_4}|$.
\begin{tikzpicture}[line width = 1pt]
   \coordinate [label=left:{$\bfx_1$}](A) at (0.0,0.0);
   \coordinate [label=right:{$\bfx_2$}](B) at (5.0,0.0);
   \coordinate [label=above:{$\bfx_4$}](D) at (1.4,3.2);
   \coordinate [label=right:{\raisebox{10pt}{$\bfx_3$}}](C) at (1.5,0.9);
   \draw (A) to node[below]{$\alpha_1$}(B);
   \draw (B) to node[above]{$\delta$}(D);
   \draw (A) to node[left]{$\alpha_3$} (D);
   \draw[dotted] (D) to node[right]{$\gamma$}(C);
   \draw[dotted] (C) to node[above]{$\beta$}(B);
   \draw[dotted] (A) to node[above]{$\alpha_2$} (C);
\end{tikzpicture}
\hspace{2mm}
\raisebox{2.1cm}{
\begin{minipage}[c]{9.4cm}
By the definition of the standard position, we have 
\begin{align*}
 \alpha_2 \le \min\{\alpha_3, \beta,\gamma\} \le 
 \max\{\alpha_3, \beta,\gamma\} \le \alpha_1.
\end{align*}
Hence, we have either $h_T = \alpha_1$ or $h_T = \delta$.
Note that $\overline{\bfx_1\bfx_4}$ is the shortest edge of the triangle
$\triangle \bfx_1\bfx_2\bfx_4$ because $\bfx_1$ and $\bfx_4$ belong to
the same half-space.
\end{minipage}
}
Hence, we have $\alpha_3 \le \delta$ and
\begin{align*}
   \alpha_1 \le h_T < 2 \alpha_1,  \quad \text{ or } \quad
   \frac{1}{2}h_T < \alpha_1 \le h_T.
\end{align*}
So far, we realize that either $h_2 = \alpha_3$, $h_2 = \beta$, or
$h_2 = \gamma$.   Recall that $\alpha_2 = h_1$.  In the following,
we check each case.
\begin{itemize}
 \item Case of $h_2 = \alpha_3$.  In this case, we have
$\alpha_1 \alpha_2 \alpha_3 = \alpha_1 h_1 h_2$, and
\begin{align*}
   \alpha_1 \alpha_2 \alpha_3 \le  h_1 h_2 h_T <
   2 \alpha_1 \alpha_2 \alpha_3  \quad \text{ amd } \quad
   H_T \le R_T < 2 H_T.
\end{align*}
\item Case of $h_2 = \beta$.  Note that
$h_2 = \beta \le \alpha_3$, and $\overline{\bfx_1\bfx_2}$ and
$\overline{\bfx_1\bfx_3}$ are the longest and shortest edges of
$\triangle \bfx_1\bfx_2\bfx_3$, respectively.  Therefore, we have
\begin{align*}
  \frac{1}{2}\alpha_3 \le \frac{1}{2}\alpha_1 
  < \beta = h_2 \le \alpha_3 \le \alpha_1.
\end{align*}
This means that
\begin{align*}
  \frac{1}{2}\alpha_1 \alpha_2 \alpha_3 < h_1h_2h_T \le 
  2 \alpha_1 \alpha_2 \alpha_3 \quad \text{ and } \quad
  \frac{1}{2}H_T < R_T \le 2H_T.
\end{align*}
\item Case of $h_2 = \gamma$.  Note that
$h_2 = \gamma \le \alpha_3$, and $\overline{\bfx_1\bfx_4}$ and
$\overline{\bfx_1\bfx_3}$ are the longest and shortest edges of
$\triangle \bfx_1\bfx_3\bfx_4$, respectively.  Therefore, we have
\begin{align*}
  \frac{1}{2}\alpha_3  < \gamma = h_2 \le \alpha_3.
\end{align*}
This implies
\begin{align*}
  \frac{1}{2}\alpha_1 \alpha_2 \alpha_3 < h_1h_2h_T \le 
  2 \alpha_1 \alpha_2 \alpha_3 \quad \text{ and } \quad
  \frac{1}{2}H_T < R_T \le 2H_T.
\end{align*}
\end{itemize}

\noindent
\textbf{Case~2.} Suppose that $T$ is of Type~2.  Set
 $\beta := |\bfx_1\bfx_3|$, $\gamma := |\bfx_3\bfx_4|$, and
$\delta := |\bfx_2\bfx_4|$.
\begin{tikzpicture}[line width = 1pt]
   \coordinate [label=left:{$\bfx_1$}](A) at (0.0,0.0);
   \coordinate [label=right:{$\bfx_2$}](B) at (5.0,0.0);
   \coordinate [label=above:{$\bfx_4$}](D) at (1.6,3.0);
   \coordinate [label=right:{\raisebox{5pt}{$\bfx_3$}}](C) at (4.0,1.3);
   \draw (A) to node[below]{$\alpha_1$}(B);
   \draw (B) to node[right]{$\alpha_2$}(C);
   \draw (C) to node[above]{$\gamma$} (D);
   \draw (D) to node[below]{$\delta$} (B);
   \draw (A) to node[left]{$\alpha_3$} (D);
   \draw[dotted] (A) to node[above]{$\beta$} (C);
\end{tikzpicture}
\hspace{2mm}
\raisebox{2cm}{
\begin{minipage}[c]{9.4cm}
By the definition of the standard position, we have 
\begin{align*}
 \alpha_2 \le \min\{\beta,\gamma,\delta\} \le 
 \max\{\beta,\gamma,\delta\} \le \alpha_1.
\end{align*}
Note that $\overline{\bfx_1\bfx_2}$ is the longest edge of the triangle
$\triangle \bfx_1\bfx_2\bfx_4$ because $\bfx_1$ and $\bfx_4$ belong to
the same half-space.  Hence, we have
$\alpha_3 \le \delta \le \alpha_1 = h_T$.
\end{minipage}
}
Therefore, we realize that either $h_2 = \alpha_3$, $h_2 = \beta$, or
$h_2 = \gamma$.   In the following, we check each case.
\begin{itemize}
 \item Case of $h_2 = \alpha_3$.  In this case, we have
$\alpha_1 \alpha_2 \alpha_3 = h_1h_2h_T$ and
$H_T = R_T$.
\item Case of $h_2 = \beta$.  Note that
$h_2 = \beta \le \alpha_3$, and $\overline{\bfx_1\bfx_2}$ and
$\overline{\bfx_2\bfx_3}$ are the longest and shortest edges of
$\triangle \bfx_1\bfx_2\bfx_3$, respectively.  Therefore, we have
\begin{align*}
  \frac{1}{2}\alpha_3 \le \frac{1}{2}\alpha_1 
  < \beta = h_2 \le \alpha_3 \le \alpha_1.
\end{align*}
This implies
\begin{align*}
  \frac{1}{2}\alpha_1 \alpha_2 \alpha_3 < h_1h_2h_T \le 
  \alpha_1 \alpha_2 \alpha_3 \quad \text{ and } \quad
  \frac{1}{2}H_T < R_T \le H_T.
\end{align*}
\item Case of $h_2 = \gamma$.  Note that
$h_2 = \gamma \le \alpha_3 \le \delta$,
 and $\overline{\bfx_2\bfx_4}$ and
$\overline{\bfx_1\bfx_3}$ are the longest and shortest edges of
$\triangle \bfx_2\bfx_3\bfx_4$, respectively.  Therefore, we have
\begin{align*}
  \frac{1}{2}\alpha_3  \le \frac{1}{2}\delta < \gamma = h_2
  \le \alpha_3 \le \delta.
\end{align*}
This implies
\begin{align*}
  \frac{1}{2}\alpha_1 \alpha_2 \alpha_3 < h_1h_2h_T \le 
   \alpha_1 \alpha_2 \alpha_3 \quad \text{ and } \quad
  \frac{1}{2}H_T < R_T \le H_T.
\end{align*}
\end{itemize}
Therefore, all cases are checked and the proof is completed. 
$\square$

\vspace{3mm}
\noindent
\begin{figure}[b]
\centering
\includegraphics[width=8cm]{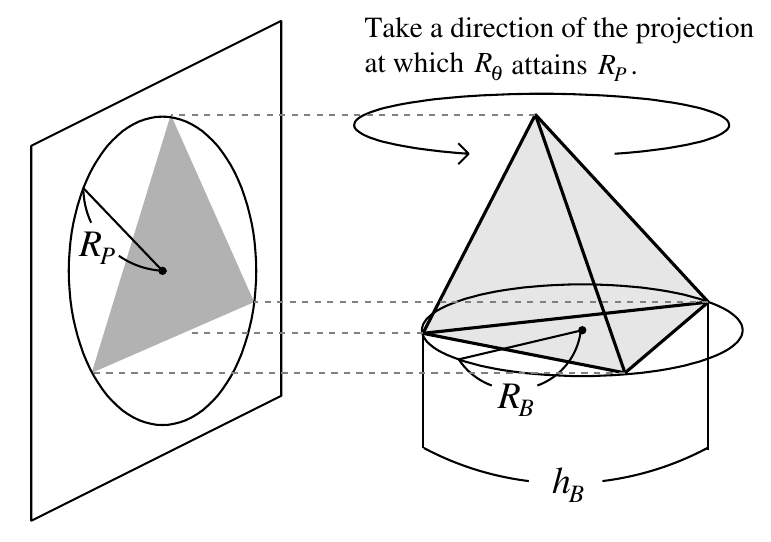}
\caption{The image of the projected circumradius of $T$.}
\label{projected-circum}
\end{figure}
\textit{Remark.} In \cite{KobayashiTsuchiya5}, the \textit{projected
circumradius} $\widetilde{R}_T$ is defined for a tetrahedron $T$ as
follows. Take any facet $B$ of $T$, and suppose that $T$ is transformed
by translation and rotation so that $B$ is on $xy$-plain.
Let $P_{xz}$ be the perpendicular projection of $\R^3$ onto
$xz$-plain; $P_{xz}(x,y,z) := (x, 0, z)$.  Note that the image
$P_{xy}(T)$ is a triangle, and let $R_0$ be its circumradius.
Now, consider rotating $T$ around the circumcenter of $B$ on
$xy$-plain. Let $T_\theta$ be the rotated tetrahedron, where $\theta$ is
the angle of the rotation.  Let $R_\theta$ be the circumradius
of $P_{xz}(T_\theta)$ (see Figure~\ref{projected-circum}).
Then, define
\begin{align}
 R_P := \max_{\theta \in [-\pi/2,\pi/2]} R_\theta, \qquad
 \widetilde{R}_T := \min_{B} \frac{R_PR_B}{h_B},
\end{align}
where $R_B$ is the circumradius of $B$, $h_B := \mathrm{diam}B$,
and the minimum is taken over all the facets of $T$.
In \cite{KobayashiTsuchiya5}, a theorem similar to
Theorem~\ref{main-thm} is proved using $\widetilde{R}_T$. It is
conjectured that $R_T$ defined by \eqref{def-RT} and the projected
circumradius $\widetilde{R}_T$ are equivalent. 

While the circumradius of a triangle is a good and simple geometric
quantity that represent its ``badness'' (or ``goodness''),
it is not so clear what is the \textit{best} geometric
quantity of a tetrahedron that represents its ``badness''.

\subsection{Squeezing theorem}
As is explained in Section~\ref{stan-posi}, we may assume without loss
of generality that an arbitrary tetrahedron $T$ may be in the standard
position.  Let $\Tabc := D\bT$, where the diagonal matrix $D$ is defined
in \eqref{diag-mat}.  We define the set
$\T_p^{k}(T) \subset W^{k+1,p}(T)$ by
\begin{gather*}
   \T_p^{k}(T) := \left\{v \in W^{k+1,p}(T) \Bigm|
     v(\bfx) = 0, \; \forall \bfx \in \Sigma^k(T) \right\}.
\end{gather*}
Then, we have the following \textbf{squeezing theorem}.
\begin{center}
\fbox{
\begin{minipage}{15truecm}
\begin{theorem}\label{squeezingtheorem}
Let $k$ and $m$ be integers with $k \ge 1$ and $0 \le m \le k$.
Let $p$ be taken as \eqref{p-cond}.
Then, there exists a constant
$C_{k,m,p}$ depending on $k$, $m$, $p$, but independent of
$\alpha_i$ $(i = 1, 2, 3)$ such that
\begin{gather*}
   B_p^{m,k}(\Tabc) := \hspace{-0.2cm}
   \sup_{v \in \T_p^{k}(\Tabc)}
  \frac{|v|_{m,p,\Tabc}}{|v|_{k+1,p,\Tabc}}
   \le \left(\max_{i=1,2,3} \alpha_i\right)^{k+1-m}C_{k,m,p}.
\end{gather*}
\end{theorem}
\end{minipage}
}
\end{center}
\noindent
\textit{Proof.} Because the proof is very similar to that of
\cite[Theorem~21]{KobayashiTsuchiya6}, we give it in 
Appendix.  $\square$

\section{Proof of Theorem~\ref{main-thm}}
In this section, we prove Theorem~\ref{main-thm} using the setting
prepared so far.  Suppose that an arbitrary tetrahedron $T$ is in the
standard position.  Recall that $T = AD(\bT)$ and $\Tabc := D\bT$,
where $(A,\bT) = (\widehat{A},\hT)$ or
 $(A,\bT) = (\widetilde{A},\tT)$ defined by \eqref{diag-mat}
according to the type of $T$.
Let $v \in W^{k+1,p}(T)$, and $\tv \in W^{k+1,m}(\Tabc)$
be defined by $\tilde{v}(x) = v(A\bfx)$.  Then, it follows from
\cite[Lemma~12]{KobayashiTsuchiya6} that
\begin{gather*}
  |v|_{m,p,T} \le 3^{m\mu(p)}  t^{1/p}
    \|A^{-1}\|^{m} |\tv|_{m,p,\Tabc}, \\
  3^{-(k+1) \mu(p)} t^{1/p} \|A\|^{-(k+1)}
   |\tv|_{k+1,p,\Tabc} \le |v|_{k+1,p,T}.
\end{gather*}
Combining the above inequalities and
Theorem~\ref{squeezingtheorem}, we obtain
\begin{align*}
   \frac{|v|_{m,p,T}}{|v|_{k+1,p,T}} & \le
   c_{k,m,p} \|A\|^{k+1}\|A^{-1}\|^{m}
   \frac{|\tv|_{m,p,\Tabc}}{|\tv|_{k+1,p,\Tabc}} \\
  & \le c_{k,m,p}C_{k,m,p} \|A\|^{k+1} \|A^{-1}\|^{m}
   \left(\max_{i=1,2,3}\alpha_i\right)^{k+1-m} \\
  & \le c_{k,m,p}C_{k,m,p} \|A\|^{k+1} \|A^{-1}\|^{m}h_T^{k+1-m}.
\end{align*}
where $c_{k,m,p} := 3^{(k+1+m)\mu(p)}$. Therefore, we obtain
the following lemma.

\begin{center}
 \fbox{
\begin{minipage}{15truecm}
\begin{lemma}  \label{liu-kikuchi}
For an arbitrary triangle $T$ in the standard position, we have 
\begin{align*}
   B_p^{m,k+1}(T) := \sup_{v \in \T_p^{k}(T)}
    \frac{|v|_{m,p,T}}{|v|_{k+1,p,T}} \le c_{k,m,p} C_{k,m.p}
  \|A\|^{k+1}\|A^{-1}\|^{m}h_T^{k+1-m}.
\end{align*} 
Therefore, inserting $v - \I_T^k v \in \T_p^k(T)$ into $v$, we have
\begin{align*}
 |v - \I_T^k v|_{m,p,T} \le c_{k,m,p} C_{k,m.p}
  \|A\|^{k+1}\|A^{-1}\|^{m} 
   h_T^{k+1-m} |v|_{k+1,p,T}, \quad
   \forall v \in W^{k+1,p}(T).
\end{align*} 
\end{lemma}
\end {minipage}
}
\end{center}

We attempt to obtain upper bounds of $\|A\|$ and $\|A^{-1}\|$.
From \eqref{mat-norm}, \eqref{mat-norm2}, \eqref{def-RT}, and
\eqref{equiv}, we know that
\begin{align*}
  \|A\| \le 2, \qquad \|A^{-1}\| \le \frac{2}{t_1t_2} =
   \frac{H_T}{3h_T} \le \frac{2R_T}{3h_T}.
\end{align*}
Hence, redefining the constant $C_{k,m,p}$ (recall that 
the Sobolev (semi-)norms may be affected by rotation up to a
constant \cite[(16)]{KobayashiTsuchiya6}),
Theorem~\ref{main-thm} is proved.

\section{Proof of Theorem~\ref{thm-max-angle}}
In this section, we give a proof of Theorem~\ref{thm-max-angle}.
For the proof, we introduce the following notation convention on $T$.
Let $F_i$ be the face of $T$ opposite to $\bfx_i$.  We denote  the dihedral
angle between the faces $F_i$ and $F_j$ by $\psi^{i,j}$. Note that
$\psi^{i,j} = \psi^{j,i}$. Furthermore, we denote the internal angle at
$\bfx_j$ on $F_i$  by $\theta_j^i$, and the angle between $F_i$ and
$\overline{\bfx_i\bfx_j}$ by $\phi_j^i$. 

\begin{center}
Table~1. Notation convention on $T$ 
$(i,j = 1, 2, 3, 4, \; i \neq j)$. \\[6pt]
\begin{tabular}{|c|l|}\hline
  $\bfx_i$  & the vertices of $T$. \\ \hline
  $F_i$  & the face opposite to $\bfx_i$. \\ \hline
  $\psi^{i,j}$ 
   &  the dihedral angle between $F_i$ and $F_j$. \\ \hline
  $\theta_j^i$ 
   &  the internal angle of $F_i$ at $\bfx_j$. \\ \hline
  $\phi_j^i$ 
   & the angle between $F_i$ and $\overline{\bfx_i \bfx_j}$. \\ \hline
 \end{tabular}
\end{center}

\begin{figure}[htbp]
\begin{tikzpicture}[line width = 1pt]
   \coordinate [label=below:{$\bfx_n$}](A) at (0.5,0.0);
   \coordinate [label=right:{$\bfx_k$}](B) at (4.5,1.5);
   \coordinate [label=above:{$\bfx_j$}](D) at (2.0,5.0);
   \coordinate [label=left:{$\bfx_m$}](C) at (-1.0,2.0);
   \draw[line width=1pt] (A) -- (B) -- (D) -- (C) -- (A) -- (D);
   \draw[dotted] (B) -- (C);
   \draw[thin] ([shift={(A)}] 20:0.7) 
    arc[radius=0.7, start angle=20, end angle=75]
    node[right]{$\quad \theta_n^m$};
   \coordinate (G) at ($(A)!0.13!(D)$);
   \coordinate (H) at ($(A)!0.2!(C)$);
   \draw[line width=0.3pt] (G) to[bend right]
    node[above]{{$\theta_n^k$}} (H);
\end{tikzpicture}
\qquad
\begin{tikzpicture}[line width = 1pt,scale=1.0]
   \coordinate [label=below:{$\bfx_n$}](A) at (0.5,0.0);
   \coordinate [label=right:{$\bfx_k$}](B) at (4.5,1.5);
   \coordinate [label=above:{$\bfx_j$}](D) at (2.0,5.0);
   \coordinate [label=left:{$\bfx_m$}](C) at (-1.0,2.0);
   \coordinate [label=below:{{$A$}}](E) at (1.5,1.3);
   \coordinate [label=below:{{$B$}}](F) at ($(A)!0.55!(B)$);
   \draw[line width=1pt] (A) -- (B) -- (D) -- (C) -- (A) -- (D);
   \draw[dotted] (B) -- (C);
   \draw[line width=0.3pt] (A) -- (E) -- (D);
   \draw[line width=0.3pt] (E) -- (F) -- (D);
   \coordinate (G) at ($(A)!0.15!(D)$);
   \coordinate (H) at ($(A)!0.3!(E)$);
   \draw[line width=0.3pt] (G) to[bend left]
    node[right]{\!\raisebox{35pt}{$\phi_n^j$}} (H);
   \coordinate (P) at ($(F)!0.2!(E)$);
   \coordinate (Q) at ($(F)!0.1!(D)$);
   \draw[line width=0.3pt] (P) to[bend left]
    node[left]{\raisebox{20pt}{$\psi^{j,m}$}\!\!} (Q);
\end{tikzpicture}
\caption{Definitions of the angles on $T$.} \label{fig2}
\end{figure}
Let $A$ and $B$ be the feet of perpendicular lines from $\bfx_j$ to $F_j$
and from $\bfx_j$ to $\overline{\bfx_n\bfx_k}$, respectively (see Figure~\ref{fig2}).
Then, we have
\begin{align*}
  |\overline{\bfx_j\bfx_n}| \sin \phi_n^j = |\overline{\bfx_jA}|
   = |\overline{\bfx_jB}| \sin \psi^{j,m}
   = |\overline{\bfx_j\bfx_n}| \sin \theta_n^m \sin \psi^{j,m}.
\end{align*}
A similar equation holds for $\phi_n^j$, $\theta_n^k$, and $\psi^{k,j}$.
Therefore,
\begin{align}
  \begin{aligned}
  &  \sin \phi_n^j = \sin \theta_n^k \sin \psi^{k,j} =
    \sin \theta_n^m \sin \psi^{m,j} \\
  & \qquad j = 1, 2, 3, 4,  \quad m, n, k
    \in \{1,2,3,4\} \backslash \{j\}.
   \end{aligned}
  \label{twosin}
\end{align}
In the following, we abbreviate ``maximum angle condition''
as MAC.
 
\begin{lemma}[Cosine rules on tetrahedrons]
Let $T \subset \mathbb{R}^3$ be a tetrahedron. 
Let $j = 1, 2, 3, 4$ and $\{k, m, n\} = \{1,2,3,4\} \backslash \{j\}$.
Then, we have
\begin{align}
\cos \theta^{k}_j & = \cos \theta^{m}_j \cos \theta^{n}_j
      + \sin \theta^{m}_j \sin \theta^{n}_j \cos \psi^{m,n},
   \notag \\
\cos \psi^{{n},{m}} & = \sin \psi^{m,k} \sin \psi^{n,k}
   \cos \theta^{k}_j - \cos \psi^{m,k} \cos \psi^{n,k}.
   \label{semi17d}
\end{align}
\end{lemma}
\noindent
\textit{Proof.}
See \cite{GelCot,Tod1886}.  $\square$

\begin{lemma} \label{lem53}
Let $T \subset \mathbb{R}^2$ be a triangle and let $\theta_i$
$(i=1,2,3)$  be the internal angles of $T$ with
$\theta_1 \leq \theta_2 \leq \theta_3$. If there exists
$\gmax \in [\pi/3, \pi)$ such that
$\theta_{3} \leq \gmax$, then we have
\begin{align}
\sin \theta_{2}, \  \sin \theta_{3} \geq \min \left\{
   \sin \frac{\pi - \gmax}{2}, \sin \gmax \right\}.
  \label{eq15}
\end{align}
\end{lemma}
\noindent
\textit{Proof.} Because $\theta_1 + \theta_2 + \theta_3 = \pi$,
the assumptions yield 
\begin{align*}
 2\theta_2 \ge \theta_1 + \theta_2
 = \pi - \theta_3 \ge \pi - \gmax \quad
  \text{ and } \quad
   \frac{\pi - \gmax}{2} \le \theta_2 \le \theta_3 \le \gmax,
\end{align*}
which implies \eqref{eq15}.   $\square$

\begin{lemma} \label{lem55}
For $\gamma \in [{\pi}/{3}, \pi)$, we have
\begin{align*}
\displaystyle
0 < \frac{\cos \gamma + 1}{\sin \frac{\gamma}{2} + 1} \leq 1.
\end{align*}
\end{lemma}
\noindent
\textit{Proof.} This lemma can be proved immediately from
\begin{align*}
   \frac{\cos\gamma + 1}{\sin \frac{\gamma}{2} + 1}
  = 2 \left(1 - \sin\frac{\gamma}{2}\right), \qquad
   \frac{\pi}{6} \le \frac{\gamma}{2} < \frac{\pi}{2}, \qquad
   \frac{1}{2} \le \sin\frac{\gamma}{2} < 1. 
  \qquad \square
\end{align*}

\begin{lemma} \label{lem56}
Let $T \subset \R^3$ be a tetrahedron.  Suppose that $T$ satisfies the
MAC with $\gmax \in [\pi/3, \pi)$.
Additionally, assume that $\theta_n^j$ is not the minimum
angle of face $F_j = \triangle P_m P_nP_k$, and
$\theta_n^j < {\pi}/{2}$, where $j = 1,2,3,4$ and
$\{m,n,k\} = \{1,2,3,4\} \backslash \{j\}$.
Then, setting $\delta$ to
\begin{align*}
\displaystyle
\sin \delta = \left( \frac{\cos \gmax + 1}
   {\sin \frac{\gmax}{2} + 1} \right)^{1/2}, \qquad
   0 <  \delta \leq \frac{\pi}{2},
\end{align*}
we have either
\begin{align}
\displaystyle
\psi^{m,j} \geq \delta, \quad \text{or} \quad
   \psi^{k,j} \geq \delta. \label{semi18}
\end{align}
\end{lemma}
\noindent
\textit{Proof.} From Lemma~\ref{lem55}, we have
\begin{align*}
\displaystyle
0 < \frac{\cos \gmax + 1}{\sin \frac{\gmax}{2} + 1} \leq 1,
\end{align*}
and we confirm that $\delta$ is well-defined.

The proof is by contradiction.  Suppose that 
\begin{align*}
0 < \psi^{m,j} < \delta \quad \text{ and } \quad 0 < \psi^{k,j} < \delta.
\end{align*}
Then, we have
$0 < \sin  \psi^{m,j} \sin \psi^{k,j} < \sin^2 \delta$ and
$1 > \cos  \psi^{m,j} \cos \psi^{k,j} > \cos^2 \delta$. 
From Lemma \ref{lem53} and the assumption, we have
\begin{align*}
\frac{\pi - \gmax}{2} \leq \theta_n^j < \frac{\pi}{2},  \qquad
 0 < \cos  \theta_n^j \le 
  \cos \left( \frac{\pi - \gmax}{2} \right)
   = \sin \frac{\gmax}{2}.
\end{align*}
Thus, we obtain
\begin{align*}
  \sin  \psi^{m,j} \sin \psi^{k,j} \cos \theta_n^j
   < \sin^2 \delta \sin \frac{\gmax}{2}.
\end{align*}
The cosine rule \eqref{semi17d} and the above inequalities yield
\begin{align*}
\cos \psi^{m,k} & = \sin \psi^{m,j}
  \sin \psi^{k,j} \cos \theta_n^{j}
   - \cos \psi^{m,j} \cos \psi^{k,j} \\
 & < \sin^2 \delta \sin \frac{\gmax}{2} - (1 - \sin^2 \delta ) \\
&= \frac{\cos \gmax + 1}{\sin \frac{\gmax}{2} + 1}
 \left( \sin \frac{\gmax}{2} + 1 \right) - 1 =  \cos
 \gmax,
\end{align*}
which contradicts the MAC: $\psi^{m,k} \leq \gmax$.  $\square$

\begin{corollary} \label{coro57}
Under the assumptions of Lemma \ref{lem56}, we have
\begin{align*}
\sin \psi^{m,j} \ge C_0, \quad \text{or}
 \quad \sin \psi^{k,j} \ge C_0, \qquad
  C_0 := \min \{ \sin\delta , \sin\gmax \}.
\end{align*}
\end{corollary}

\begin{lemma} \label{lem59}
For $j = 1,2,3, 4$, let $\{m,n,k\} = \{1,2,3,4\} \backslash \{j\}$.
Let $p \in \{m, n, k\}$, and $\{q, r\} = \{m, n, k\} \backslash \{p\}$.
Suppose that there exists a positive constant $M$ 
with $0 < M < 1$ such that $\sin \phi_p^j \sin\theta_n^j\ge M$.
Then, setting $\gamma (M) := \pi - \sin^{-1} M$
$(\frac{\pi}{2} < \gamma(M) < \pi)$, the MAC with
$\gamma(M)$ is satisfied on faces $F_j$, $F_q$, $F_r$, and
$\psi^{j,q}$, $\psi^{j,r} \le \gamma(M)$.
\end{lemma}
\noindent
\textit{Proof.} From the assumption, we have
\begin{align*}
 M \le \sin \phi_p^j \sin\theta_n^j \le \sin \theta_n^j
   \quad \text{ and } \quad
  M \le \sin \phi_p^j.
\end{align*}
Hence, the definition of $\gamma(M)$ yields
$\pi - \gamma(M) \le \theta_n^j \le \gamma(M)$.
Because $\theta_n^j + \theta_m^j + \theta_k^j = \pi$,
we see that $\theta_m^j$,
 $\theta_k^j < \theta_m^j + \theta_k^j  \le \gamma(M)$.
That is, the MAC with $\gamma(M)$ is satisfied on
face $F_j = \triangle P_mP_nP_k$.

Moreover, it follows from \eqref{twosin} that 
\begin{align*}
  M & \le \sin \phi_p^j = \sin \theta_p^q \sin \psi^{q,j} 
      = \sin \theta_p^r \sin \psi^{r,j} \\
     & \le \sin \theta_p^q, \; \sin \theta_p^r, \;
        \sin \psi^{r,j}, \; \sin \psi^{q,j} 
\end{align*}
By the same reasoning, we find that the MAC
with $\gamma(M)$ is satisfied on faces $F_q$ and $F_r$, and
$\psi^{j,q}$, $\psi^{j,r} \le \gamma(M)$.  $\square$

\vspace{3mm}
In the following, we prove Theorem~\ref{thm-max-angle} using
$H_T$ instead of $R_T$.  We divide the proof into four cases.

\subsection{Type~1: Proof of ``MAC implies \eqref{equiv-cond}''}
\label{Case1}
First, we suppose that $T$ is of Type~1 and satisfies the MAC
with $\gmax$, $\pi/3 \le \gmax < \pi$.
Because $|T| = \frac{1}{6} \alpha_1 \alpha_2 \alpha_3
\sin \theta_1^4  \sin \phi_1^4$, we have
 \begin{align*}
   \frac{H_T}{h_T} = \frac{\alpha_1 \alpha_2 \alpha_3}{|T|}
     = \frac{6}{\sin \theta_1^4  \sin \phi_1^4}.
\end{align*}
From the definition of Type~1, we realize that
$\theta_2^4 \le \theta_1^4 \le \theta_3^4$, that is,
$\theta_3^4$ and $\theta_2^4$ are the maximum and minimum angles
of face $F_4 = \triangle P_1 P_2 P_3$, respectively. 
Thus, it follows from Lemma \ref{lem53} that
\begin{align*}
 \frac{\pi - \gmax}{2} \leq \theta_1^4 \leq \gmax, \quad  \sin
 \theta_1^4 \geq \min \left \{ \sin \frac{\pi - \gmax}{2}, \sin
 \gmax \right \} =: C_1.
\end{align*}
Additionally, we may apply Lemma~\ref{lem56} to $\theta_1^4$ and $F_4$,
and find that either $\psi^{2,4} \ge \delta$ or
$\psi^{3,4} \ge \delta$, where
$\delta = \delta(\gmax)$, $0 < \delta \le \pi/2$ is defined as
\begin{align}
  \sin \delta = \left( \frac{\cos \gmax + 1}
   {\sin \frac{\gmax}{2} + 1} \right)^{1/2}.
   \label{def-delta}
\end{align}

Suppose that $\psi^{2,4} \ge \delta$. By Corollary~\ref{coro57}
and \eqref{twosin}, we have 
\begin{align*}
  \sin \phi_1^4 = \sin \theta^2_1 \sin \psi^{2,4}
  \ge C_0 \sin \theta^2_1,
\end{align*}
where $C_0$ is the constant defined in Corollary~\ref{coro57}.
By the definition of Type 1, $\theta_1^2$ is not
the minimum angle of $F_2 = \triangle P_1 P_3 P_4$, and therefore, we
have
\begin{align*}
  \frac{\pi - \gmax}{2} \le \theta_1^2 \leq \gmax,
  \quad  \sin \theta_1^2 \ge  C_1.
\end{align*}
Thus, we obtain $\sin \phi_1^4 \geq C_0 C_1$.

Next, suppose that $\psi^{3,4} \geq \delta$.
Replacing $\psi^{2,4}$, $\theta_1^2$, and $F_2$ with
$\psi^{3,4}$, $\theta_1^3$, and $F_3$ in the above argument,
we obtain $\sin \phi_1^4 \geq C_0 C_1$ in the same manner.

Gathering the above results, we conclude that
\begin{align*}
  \frac{H_T}{h_T} = \frac{6}{\sin \theta_1^4  \sin \phi_1^4}
  \le \frac{6}{C_0 C_1^2} =: D
\end{align*}
in both cases, that is, \eqref{equiv-cond} holds.

\subsection{Type~1: Proof of ``\eqref{equiv-cond} implies MAC''}
\label{Case2}
Now, we suppose that $T$ is of Type~1 and 
\begin{align*}
  \frac{H_T}{h_T} = \frac{\alpha_1 \alpha_2 \alpha_3}{|T|}
   = \frac{6}{\sin \theta_1^4  \sin \phi_1^4}  \le D.
\end{align*}
Because $\theta^4_1 < \pi/2$ and 
$\sin \theta_1^4  \sin \phi_1^4 < 1$, we have
\begin{align*}
\sin \theta_1^4  \sin \phi_1^4 \geq \frac{6}{D} =: M, \qquad
  0 < M < 1.
\end{align*}
By Lemma~\ref{lem59} with $j = 4$ and $p = 1$, setting 
$\gamma (M) := \pi - \sin^{-1} M$, we have 
$\frac{\pi}{2} < \gamma(M) < \pi$, and the MAC
with $\gamma(M)$ is satisfied on $F_2$, $F_3$, $F_4$, and 
$\psi^{2,4}$, $\psi^{3,4} \le \gamma(M)$.

Note that $|T| = \frac{1}{6}\alpha_1\alpha_2\alpha_3
\sin \theta_1^3 \sin \phi_1^3$, and we have
\begin{align*}
  \frac{H_T}{h_T} = \frac{\alpha_1 \alpha_2 \alpha_3}{|T|}
   = \frac{6}{\sin \theta_1^3 \sin \phi_1^3}  \le D.
\end{align*}
Thus, by Lemma~\ref{lem59} with $j = 3$ and $p = 1$,
we find that $\psi^{2,3} \le \gamma(M)$.

Because $|\overline{P_3 P_4}| < |\overline{P_1 P_4}| +
|\overline{P_1P_3}| \leq 2 \alpha_3$
on $F_2 = \triangle P_1 P_3 P_4$
and $|\overline{P_2P_3}| \leq \alpha_1$, we note that
\begin{align*}
  |T| = \frac{1}{6} \alpha_2  |\overline{P_2 P_3}| |\overline{P_3 P_4}|
    \sin \theta^1_3 \sin \phi_3^1 
   < \frac{1}{3} \alpha_1 \alpha_2 \alpha_3 \sin \theta^1_3 \sin \phi_3^1.
\end{align*}
Thus, we have
\begin{align*}
  D \ge \frac{H_K}{h_K} > \frac{3}{\sin \theta^1_3 \sin \phi_3^1}
  \quad \text{ and } \quad
  \sin \theta^1_3 \sin \phi_3^1 > \frac{3}{D} = \frac{M}{2}.
\end{align*}
From Lemma~\ref{lem59}, setting $\gamma (M/2) := \pi - \sin^{-1} (M/2)$,
we have $\frac{\pi}{2} < \gamma(M/2) < \pi$ and 
MAC with $\gamma(M/2)$ is satisfied on $F_1$,
and $\psi^{2,1}$, $\psi^{4,1} \le \gamma(M/2)$.

The final thing to prove is the MAC for $\psi^{1,3}$.
From the cosine rule \eqref{semi17d}, we have
\begin{align*}
\cos \psi^{1,3} = \sin \psi^{3,4} \sin \psi^{4,1}
   \cos \theta_2^{4} - \cos \psi^{3,4} \cos \psi^{4,1}.
\end{align*}
By the definition of Type~1, the angle $\theta_2^4$ is the minimum angle
of $F_4 = \triangle P_1 P_2 P_3$, and therefore, we have
\begin{align*}
   \cos \theta_2^4 \ge \frac{1}{2}, \;
\sin \psi^{{3},{4}} \sin \psi^{{4},{1}} \cos \theta^{4}_2 > 0,
\quad\text{ and }
 \cos \psi^{1,3} > - \cos \psi^{3,4} \cos \psi^{4,1}.
\end{align*}
From the above argument, we have $\sin \psi^{3,4} > M$,
$\sin \psi^{4,1} > M/2$, and
\begin{align*}
\cos \psi^{1,3} & > - \cos \psi^{3,4} \cos \psi^{4,1} 
 \ge - | \cos \psi^{3,4} | | \cos \psi^{4,1} | \\
 & = - \sqrt{1 - \sin^2 \psi^{3,4}} \sqrt{1 - \sin^2 \psi^{4,1}}
 > - \sqrt{1 - M^2} \sqrt{1 - \frac{M^2}{4}} > -1.
\end{align*}
Therefore, we conclude that
\begin{align*}
\psi^{1,3} < \cos^{-1}
  \left(- \sqrt{1 - M^2} \sqrt{1 - \frac{M^2}{4}}\right) < \pi,
\end{align*}
and $T$ satisfies the MAC with
\begin{align*}
  \gmax := \max\left\{\gamma(M/2),  \cos^{-1}
  \left(- \sqrt{1 - M^2} \sqrt{1 - \frac{M^2}{4}}\right)
   \right\}.
\end{align*}

\subsection{Type~2: Proof of ``MAC implies \eqref{equiv-cond}''}
First, we suppose that $T$ is of Type~2 and satisfies the MAC
with $\gmax \in [\pi/3, \pi)$.  The proof is very similar to
that described in Section~\ref{Case1}.

By the definition of Type~2,
$\alpha_3 = |\overline{P_1P_4}| \le |\overline{P_2P_4}|$.
Because
\begin{align*}
    |T| = \frac{1}{6} \alpha_1 \alpha_2 \alpha_3
   \sin \theta_2^4  \sin \phi_1^4 =
    \frac{1}{6} \alpha_1 \alpha_2 |\overline{P_2P_4}|
    \sin \theta_2^4  \sin \phi_2^4,
\end{align*}
we have
 \begin{align}
   \frac{H_T}{h_T} = \frac{\alpha_1 \alpha_2 \alpha_3}{|T|}
     = \frac{6}{\sin \theta_2^4  \sin \phi_1^4}
  \le \frac{6}{\sin \theta_2^4  \sin \phi_2^4}.
  \label{Case3}
\end{align}
From the definition of Type~2, we realize that
$\theta_1^4 \le \theta_2^4 \le \theta_3^4$ on $F_4$, 
$\theta_2^3 \le \theta_1^3 \le \theta_4^3$ on $F_3$, and
$\theta_2^1$ is not the minimum angle of $F_1$.
Thus, it follows from Lemma \ref{lem53} that
\begin{align*}
 \frac{\pi - \gmax}{2} \leq \theta_2^4, \, \theta_1^3, 
  \, \theta_2^1\le \gmax,
  \quad  \sin \theta_2^4, \, \sin \theta_1^3, \,
  \sin \theta_2^1 \ge C_1.
\end{align*}
Additionally, we may apply Lemma~\ref{lem56} to $\theta_2^4$ and $F_4$,
and find that either $\psi^{1,4} \ge \delta$ or
$\psi^{3,4} \ge \delta$, where $\delta = \delta(\gmax)$ is defined by
\eqref{def-delta}.

Suppose that $\psi^{3,4} \ge \delta$.  Using the same argument as in
Section~\ref{Case1}, we have
\begin{align*}
  \sin \phi_1^4 = \sin \theta_1^3 \sin \psi^{3,4}
  \ge C_0 \sin \theta_1^3 \ge C_0C_1.
\end{align*}

Next, suppose that $\psi^{1,4} \ge \delta$.  We have 
\begin{align*}
  \sin \phi_2^4 = \sin \theta_2^1 \sin \psi^{1,4}
  \ge C_0 \sin \theta_2^1 \ge C_0C_1.
\end{align*}

Combining these results with \eqref{Case3}, we obtain
\begin{align*}
  \frac{H_T}{h_T}  \le \frac{6}{C_0 C_1^2} =: D,
\end{align*}
that is, \eqref{equiv-cond} holds.

\subsection{Type~2: Proof of ``\eqref{equiv-cond} implies MAC''}
Finally, we suppose that $T$ is of Type~2 and 
\begin{align*}
  \frac{H_T}{h_T} = \frac{\alpha_1 \alpha_2 \alpha_3}{|T|}
   = \frac{6}{\sin \theta_2^4  \sin \phi_1^4} \le D, \quad
   \sin \theta_2^4  \sin \phi_1^4 \ge \frac{6}{D} =: M. 
\end{align*}
The proof is very similar to that described in Section~\ref{Case2}.
By Lemma~\ref{lem59} with $j = 4$ and $p = 1$, setting 
$\gamma (M) := \pi - \sin^{-1} M$, 
the MAC with $\gamma(M)$ is satisfied on
$F_2$, $F_3$, $F_4$, and 
$\psi^{2,4}$, $\psi^{3,4} \le \gamma(M)$.

Because $|\overline{P_2 P_4}| \le \alpha_1$, we have
\begin{align*}
  |T| = \frac{1}{6} |\overline{P_2 P_3}|
  |\overline{P_2 P_4}| |\overline{P_1 P_4}|
    \sin \theta_2^1 \sin \phi_4^1 
   \le \frac{1}{6} \alpha_1 \alpha_2 \alpha_3 \sin \theta_2^1
    \sin \phi_4^1.
\end{align*}
This yields
\begin{align*}
  D \ge \frac{H_T}{h_T} \ge \frac{6}{\sin \theta_2^1 \sin \phi_4^1}
  \quad \text{ and } \quad
  \sin \theta_2^1 \sin \phi_4^1 \ge \frac{6}{D} = M,
\end{align*}
and, by Lemma~\ref{lem59} with $j = 1$ and $p = 4$, we find 
that the MAC with $\gamma(M)$ is satisfied on $F_1$,
and $\psi^{1,2}$, $\psi^{1,3} \le \gamma(M)$.

The final thing to prove is the MAC for
$\psi^{1,4}$ and $\psi^{2,3}$.  
By the cosine rule \eqref{semi17d} with $j=2$, we have
\begin{align*}
  \cos \psi^{1,4} & = \sin \psi^{1,3} \sin \psi^{4,3}
   \cos \theta_2^3 - \cos \psi^{1,3} \cos \psi^{4,3}, \\
   \cos \psi^{2,3} & = \sin \psi^{2,4} \sin \psi^{3,4}
   \cos \theta_1^4 - \cos \psi^{2,4} \cos \psi^{3,4}.
\end{align*}
By the definition of Type 2,
$\theta_2^3$ and $\theta_1^4$ are the minimum angles of $F_3$
and $F_4$, respectively.  Therefore, we have
$\cos \theta_2^3$, $\cos \theta_1^4 \ge \frac{1}{3}$ and thus
\begin{align*}
 \cos \psi^{1,4} > - \cos \psi^{1,3} \cos \psi^{3,4}, \qquad
 \cos \psi^{2,3} > - \cos \psi^{2,4} \cos \psi^{3,4}.
\end{align*}
Because $\sin \psi^{1,3}$, $\sin \psi^{2,4}$, 
$\sin \psi^{3,4} > M$, we find that
\begin{align*}
   \cos \psi^{1,4} & > - \cos \psi^{1,3} \cos \psi^{3,4} 
    \ge - \sqrt{1 - \sin^2 \psi^{1,3}} \sqrt{1 - \sin^2 \psi^{3,4}}
    > M^2 - 1, \\
   \cos \psi^{2,3} & > M^2 - 1.
\end{align*}
Therefore, we conclude that
$\psi^{1,4}$, $\psi^{2,3} < \cos^{-1}(M^2 - 1) < \pi$,
and $T$ satisfies the MAC with
\begin{align*}
  \gmax := \max\left\{\gamma(M), \cos^{-1}(M^2 - 1)
   \right\}.
\end{align*}


\section*{Appendix: Proof of Theorem~\ref{squeezingtheorem}}
The proof of Theorem~\ref{squeezingtheorem} is very similar to
that of \cite[Theorem~13]{KobayashiTsuchiya4} and
\cite[Theorem~21]{KobayashiTsuchiya6}.
First, refer to \cite[Section~5]{KobayashiTsuchiya6} for
the definition of difference quotients of one and two variable
functions.  Difference quotients of three variable functions is their
simple extension. 

For a positive integer $k$, $X^k$ is the set of lattice points
defined by
\begin{align*}
   X^k & := \left\{ \bfx_{\gamma} := \frac{\gamma}{k}
     \in \R^3 \biggm|  \gamma \in \N_0^3 \right\},
\end{align*}
where $\gamma/k=(a_1/k,a_2/k,a_3/k)$ is understood as the
coordinate of a point in $\R^3$.
For $\bfx_{\gamma} \in X^k$ and a multi-index $\delta \in \N_0^3$,
we define the correspondence $\Delta^\delta$
between nodes by
$\Delta^\delta \bfx_{\gamma} := \bfx_{\gamma+\delta} = (\gamma+\delta)/k$.

For two multi-indexes  $\eta=(m_1,m_2,m_3)$,
$\delta=(n_1,n_2,n_3)$,  $\eta \le \delta$ means that
$m_i \le n_i$ $(i=1,2,3)$.
Also, $\delta\cdot\eta$ and $\delta !$ are defined by
$\delta\cdot\eta := \sum_{i=1}^3 m_i n_i$ and
 $\delta ! := n_1! n_2! n_3!$, respectively.  Suppose that, for
$\gamma, \delta \in \N_0^3$, both $\bfx_\gamma$ and
$\Delta^\delta \bfx_\gamma$ belong to $\bK$.  Then, we define the
\textit{difference quotients} for $f \in C^0(\bK)$ by
\begin{align*}
 f^{|\delta|}[\bfx_{\gamma},\Delta^\delta \bfx_{\gamma}] :=
  k^{|\delta|}\sum_{\eta \le \delta} 
\frac{(-1)^{|\delta|-|\eta|}}{\eta!(\delta-\eta)!}
    f(\Delta^\eta \bfx_{\gamma}).
\end{align*}
For example, we see that
\begin{align*}
  f^4[\bfx_{(0,0,0)},\Delta^{(2,1,1)}\bfx_{(0,0,0)}] & =
   \frac{k^4}{2} (f(\bfx_{(2,1,1)}) - 2 f(\bfx_{(1,1,1)})
   + f(\bfx_{(0,1,1)}) \\
    & \hspace{0.5cm} 
   - f(\bfx_{(2,0,1)}) + 2 f(\bfx_{(1,0,1)}) - f(\bfx_{(0,0,1)}) \\
    & \hspace{0.5cm} 
   - f(\bfx_{(2,1,0)}) + 2 f(\bfx_{(1,1,0)}) - f(\bfx_{(0,1,0)}) \\
 & \hspace{0.5cm}
   + f(\bfx_{(2,0,0)}) - 2 f(\bfx_{(1,0,0)}) + f(\bfx_{(0,0,0)})).
\end{align*}

As explained in \cite[Section~5]{KobayashiTsuchiya6}, a differential
quotients is expressed concisely by an  integral.  For that purpose, we
introduce the $s$-simplex
\begin{gather*}
   \Simp_s := \left\{(x_1,\cdots,x_s)^\top \in \R^s \mid 
   x_i \ge 0, \ 0 \le x_1 + \cdots + x_s \le 1 \right\},
\end{gather*}
and the integral of $g \in L^1(\Simp_s)$ on $\Simp_s$ is defined by
\begin{gather*}
  \int_{\Simp_s} g(w_1,\cdots,w_k) \dd\mathbf{W_s}
  := \int_{0}^{1}\int_0^{w_1}\cdots \int_0^{w_{s-1}} 
  g(w_1, \cdots, w_s)  \dd w_s \cdots \dd w_2 \dd w_1,
\end{gather*}
where $\dd \mathbf{W_s} := \dd w_1 \cdots \dd w_s$.
Then, $f^{s}[\bfx_{(l,q)},\Delta^{(0,s,0)} \bfx_{(l,q)}]$ becomes
\begin{align*}
f^{s}[\bfx_{(l,q,r)},\Delta^{(0,s,0)} \bfx_{(l,q,r)}]
  & = \int_{\Simp_s} \partial^{(0,s,0)}
  f\left(\frac{l}{k},\frac{q}{k} + \frac{1}{k}(w_1 + \cdots + w_s),
    \frac{r}{k} \right) \dd \mathbf{W_s}.
\end{align*}
For a general multi-index $(t,s,m)$, we can write
\begin{gather*}
 f^{t+s+m} [\bfx_{(l,q,r)}, \Delta^{(t,s,m)} \bfx_{(l,q,r)}] 
  = \int_{\Simp_s}\int_{\Simp_t} \int_{\Simp_m}\partial^{(t,s,m)}
  f\left(\mathbf{Z_t}, \mathbf{W_s}, \mathbf{Y_m}\right)
    \dd \mathbf{Z_t} \dd \mathbf{W_s} \dd \mathbf{Y_m}, \\
  \mathbf{Z_t} := \frac{l}{k} + \frac{1}{k}(z_1 + \cdots + z_t), \;
   \dd \mathbf{Z_t} := \dd z_1 \cdots \dd z_t, \quad
  \mathbf{W_s} := \frac{q}{k} + \frac{1}{k}(w_1 + \cdots + w_s),
     \\
  \mathbf{Y_m} := \frac{r}{k} + \frac{1}{k}(y_1 + \cdots + y_m), \quad
    \dd \mathbf{Y_m} := \dd y_1 \cdots \dd y_m.
\end{gather*}

Let $\square_{\gamma}^\delta$ be the rectangular parallelepiped defined by 
$\bfx_{\gamma}$ and $\Delta^\delta \bfx_{\gamma}$ as the diagonal
points. If $\delta=(t,s,0)$ or $(0,s,0)$,  $\square_{\gamma}^\delta$
degenerates to a rectangle or a segment.  For $v \in L^1(\hK)$
and $\square_{\gamma}^\delta$ with $\gamma=(l,q,r)$, we denote the integral as
\begin{equation*}
   \int_{\square_{\gamma}^{(t,s,m)}} v :=
   \int_{\Simp_s}\int_{\Simp_t} \int_{\Simp_m}
  v\left(\mathbf{Z_t}, \mathbf{W_s}, \mathbf{Y_m}
   \right) \dd \mathbf{Z_t} \dd \mathbf{W_s} \dd \mathbf{Y_m}.
\end{equation*}
If $\square_{\gamma}^\delta$ degenerates to a rectangle or a segment, the integral is
understood as an integral on the rectangle or on the segment.  By this notation, the
difference quotient $f^{|\delta|}[\bfx_{\gamma},\Delta^{\delta} \bfx_{\gamma}]$ is
written as 
\begin{align*}
 f^{|\delta|}[\bfx_{\gamma},\Delta^{\delta} \bfx_{\gamma}] 
  = \int_{\square_{\gamma}^{\delta}} \partial^{\delta} f.
\end{align*}
Therefore, if $u \in \T_p^k(\bT)$, then we have
\begin{align}
 0 = u^{|\delta|}[\bfx_{\gamma},\Delta^{\delta} \bfx_{\gamma}] 
  = \int_{\square_{\gamma}^{\delta}} \partial^{\delta} u, \qquad
  \forall \square_{\gamma}^{\delta} \subset \bT.
  \label{tomoko}
\end{align}

Let $S \subset \bT$ be a segment.  In the proof of
Theorem~\ref{squeezingtheorem}, the continuity of the trace operator $t$
defined as $t:W^{1,p}(\bT) \ni v \mapsto v|_S \in L^1(S)$ is crucial.
For two-dimensional case, the continuity of $t$ is standard and is
mentioned in many textbooks such as \cite{Brezis}. For three dimensional
case, the situation becomes a bit more complicated.  If the continuous
inclusion $W^{k+1,p}(\bT) \subset C^0(\bT)$ holds, the continuity of $t$
is obvious. Even if this is not the case, we still have the following
lemma.  For the proof, see \cite[Theorem~4.12]{AdamsFournier},
\cite[Lemma~2.2]{Duran}, and \cite[Theorem~2.1]{LSU}.

\begin{lemma}\label{lem31}
 Let $S \subset \bT$ be an arbitrary segment.  Then,
the following trace operators are well-defined and continuous:
\begin{align*}
   t:W^{1,p}(\bT) \to L^p(S), \quad 2 < p < \infty, \qquad
   t:W^{2,p}(\bT) \to L^p(S), \quad 1 \le p < \infty.
\end{align*}
\end{lemma}

Let $p$ be taken as \eqref{p-cond}.  The set
$\Xi_p^{\delta,k} \subset W^{k+1-|\delta|,p}(\bT)$ is then defined by
\begin{align*}
  \Xi_p^{\delta,k} & := \left\{ v \in W^{k+1-|\delta|,p}(\bT) \Bigm| 
    \int_{\square_{\gamma}^{\delta}} v = 0,\quad \forall
   \square_{\gamma}^{\delta} \subset \bT \right\}.
\end{align*}
Note that $u \in \T_p^{k}(\bT)$ implies
$\partial^\delta u \in \Xi_p^{\delta,k}$ by \eqref{tomoko}.

\begin{figure}[htbp]
\centering
{\includegraphics[width=4.05cm]{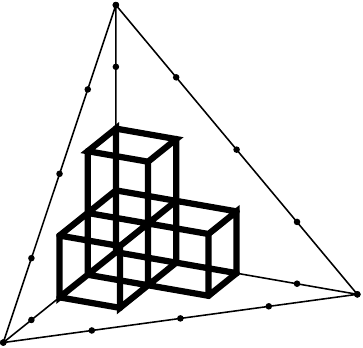}} \qquad
{\includegraphics[width=2.7cm]{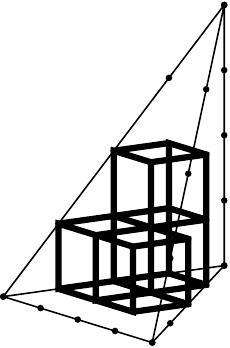}} \qquad
{\includegraphics[width=2.7cm]{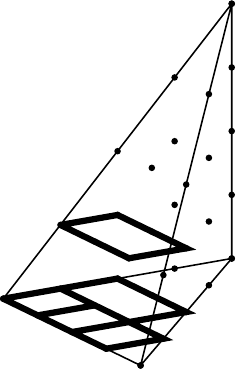}}
\caption{The four cubes and four rectangles in $\bT$.}
\label{cubes_tetra}
\end{figure}

\begin{lemma}\label{lem32} 
We have $\Xi_p^{\delta,k} \cap \mathcal{P}_{k-|\delta|} = \{0\}$.
That is, if $q \in \mathcal{P}_{k-|\delta|}$ belongs to $\Xi_p^{\delta,k}$,
then $q=0$.
\end{lemma}
\textit{Proof.}
 Note that
$\mathrm{dim} \mathcal{P}_{k-|\delta|} = \#\{\square_{\gamma}^\delta \subset \bT \}$.
For example, if $k=4$ and $|\delta|=3$, then
$\mathrm{dim}\mathcal{P}_1 = 4$.  This corresponds to the fact that,
in $\bT$, there are four cubes of size $1/4$  for $\delta=(1,1,1)$ and
there are four rectangles of size $1/2 \times 1/4$ for $\delta=(1,2,0)$.
All their vertices (corners) belong to $\Sigma^4(\bT)$
(see Figure~\ref{cubes_tetra}). Now, suppose that
$q \in  \mathcal{P}_{k-|\delta|}$ satisfies
$\int_{\square_{\gamma}^{\delta}} q = 0$ for all
$\square_{\gamma}^{\delta} \subset \bT$.  These conditions are linearly
independent and determine $q = 0$ uniquely (see Exercise below).
\hfill $\square$

\noindent
\textbf{Exercise:} Show that the condition 
``$\int_{\square_{\gamma}^{\delta}} q = 0$ for all
$\square_{\gamma}^{\delta} \subset \bT$'' implies $q = 0$
for $q \in \PP_{k-|\delta|}$.
(Hint: (1) First, consider the case $d = 1$. For example, show the
following: if a polynomial $p \in \PP_k$ satisfies
$\int_{n}^{n+1} p(x) \dd x = 0$, $n = 0, \cdots, k$, then $p = 0$.) \\
(2) Reduce the proof of the case $d > 1$ to that of the case $d - 1$.

\vspace{0.3cm}
The constant $A_p^{\delta,k}$ is defined by
\begin{align*}
   A_p^{\delta,k} := \sup_{v \in \Xi_p^{\delta,k}} \frac{|v|_{0,p,\bT}}
  {|v|_{k+1-|\delta|,p,\bT}}.
\end{align*}
The following lemma is an extension of \cite[Lemma~2.1]{BabuskaAziz}.

\begin{lemma}\label{lem33}
Let $p$ be such that $2 < p \le \infty$ if $k+1-|\delta|=1$ or
$1 \le p \le \infty$ if $k+1-|\delta|\ge 2$.  We then have
$A_p^{\delta,k} < \infty$.
\end{lemma}
\noindent
\textit{Proof.}  The proof is by contradiction.  Suppose that
$A_p^{\delta,k} = \infty$.  Then there exists a sequence 
$\{w_k\}_{i=1}^\infty \subset \Xi_p^{\delta,k}$ such that
$|w_n|_{0,p,\bT} = 1$ and 
$\lim_{n \to \infty} |w_n|_{k+1-|\delta|,p,\bT} = 0$.
By the Bramble--Hilbert lemma \cite[Theorem~14]{KobayashiTsuchiya6},
there exists $\{q_n\}\subset \PP_{k-|\delta|}$ such that
\begin{align*}
\|w_n + q_n\|_{k+1-|\delta|,p,\bT} 
   \le \inf_{q \in \PP_{k-|\delta|}}
  \|w_n + q\|_{k+1-|\delta|,p,\bT} + \frac{1}{n}
   \le C |w_n|_{k+1-|\delta|,p,\bT} + \frac{1}{n}
\end{align*}
and $\displaystyle\lim_{n \to \infty} \|w_n + q_n\|_{k+1-|\delta|,p,\bT} = 0$.
Because $\{w_n\}\subset W^{k+1-|\delta|,p}(\bT)$ is bounded,
$\{q_n\} \subset \PP_{k-|\delta|}$ is bounded as well.
Hence, there exists a subsequence $\{q_{n_i}\}$ such that
$q_{n_i}$ converges to $\bar{q} \in \PP_{k-|\delta|}$ and
$\lim_{n_i \to \infty} \|w_{n_i} + \bar{q}\|_{k+1-|\delta|,p,\bT} = 0$.
If $\square_{lp}^\delta$ is not degenerate to a rectangle or
a segment, we have
\begin{align}
  \left|\int_{\square_{lp}^\delta} (w_{n_i} + \bar{q}) \right|
  \le \int_{\square_{lp}^\delta} |w_{n_i} + \bar{q}|
  \le C \|w_{n_i} + \bar{q}\|_{k+1-|\delta|,p,\bT} \to 0 \quad
 \text{as} \; h \to 0.
  \label{ayumi}
\end{align}
If $\square_{lp}^\delta$ is degenerate to a rectangle or
a segment, \eqref{ayumi} holds as well by Lemma~\ref{lem31}.
Because $\int_{\square_{lp}^\delta} w_{n_i} = 0$ by the definition,
we have 
\begin{align}
  0 = \lim_{n_i\to\infty} \int_{\square_{lp}^\delta} 
   (w_{n_i} + \bar{q}) = \int_{\square_{lp}^\delta} \bar{q}, \qquad
    \forall \square_{lp}^\delta \subset \bT.
\end{align}
Therefore, it follows from Lemma~\ref{lem31} that $\bar{q} = 0$.
This implies that
\[
  0 = \lim_{n_i \to \infty}\|w_{n_i}\|_{k+1-|\delta|,p,\bT} 
  \ge \lim_{n_i \to \infty}|w_{n_i}|_{0,p,\bT} = 1,
\]
which is a contradiction.  $\square$

\vspace{2mm}
Define the linear transformation
by, for $(x,y,z)^\top \in \R^3$, 
\begin{align*}
   (x^*, y^*,z^*)^\top = D_{\alpha_1\alpha_2\alpha_3} (x, y,z)^\top =
    (\alpha_1 x, \alpha_2 y, \alpha_3 z)^\top, \qquad
     \alpha_i > 0, \; i = 1, 2, 3,
\end{align*}
which the diagonal matrix $D_{\alpha_1\alpha_2\alpha_3}$ is defined by
\eqref{diag-mat}.  This linear transformation  squeezes the reference
tetrahedron $\bT$ perpendicularly to
$T_{\alpha_1\alpha_2\alpha_3} = D_{\alpha_1\alpha_2\alpha_3}\bT$.  Take an
arbitrary $v \in \T_p^k(T_{\alpha_1\alpha_2\alpha_3})$
and define $u \in \T_p^k(\bT)$ by
$u(x,y,z):= v(D_{\alpha_1\alpha_2\alpha_3}(x,y,z)^\top)$.
Let $p$ be taken as \eqref{p-cond} with $m=|\delta|$.  To make formula
concise, we introduce the following notation. 
For a multi-index $\gamma = (a, b, c) \in \N_0^3$ and a real
$t \neq 0$, and $(\alpha) := (\alpha_1,\alpha_2,\alpha_3)$,
$(\alpha)^{\gamma t} :=
 \alpha_1^{a t}\alpha_2^{b t}\alpha_3^{c t}$.
Because $u \in \T_p^k(\bT)$ and $\partial^\delta u \in \Xi_p^{\delta,k}$,
we may apply Lemma~\ref{lem33} as follows.  For $p$, $1 \le p < \infty$,
we have
{\allowdisplaybreaks
\begin{align}
  \frac{|v|_{m,p,T_{\alpha_1\alpha_2\alpha_3}}^p}
    {|v|_{k+1,p,T_{\alpha_1\alpha_2\alpha_3}}^p}
   &  = \frac{\sum_{|\gamma| = m} \frac{m!}{\gamma!}
      (\alpha)^{-\gamma p}\left|\partial^{\gamma}u\right|_{0,p,\bT}^p}
      { \sum_{|\delta|= k+1} \frac{(k+1)!}{\delta!}
     (\alpha)^{- \delta p} \left|\partial^{\delta} u  \right|_{0,p,\bT}^p}
      \notag \\
  & = \frac{ \sum_{|\gamma| = m}
    \frac{m!}{\gamma!}(\alpha)^{-\gamma p}
         \left|\partial^{\gamma}u \right|_{0,p,\bT}^p }
       { \sum_{|\gamma| = m} \frac{m!}{\gamma!}
      (\alpha)^{-\gamma p}
       \left(\sum_{|\eta| = k+1-m} \frac{(k+1-m)!}
             {\eta!(\alpha)^{\eta p} }
         \left|\partial^{\eta}
                 (\partial^{\gamma}u)\right|_{0,p,\bT}^p\right)} \notag \\
  & \le \frac{ \left(\max_{i=1,2,3}\alpha_i\right)^{(k+1-m)p}
    \sum_{|\gamma| = m}
    \frac{m!}{\gamma!}(\alpha)^{- \gamma p}
         \left|\partial^{\gamma}u \right|_{0,p,\bT}^p }
       {\sum_{|\gamma| = m} \frac{m!}{\gamma!}
      (\alpha)^{-\gamma p}
       \left(\sum_{|\eta| = k+1-m} \frac{(k+1-m)!}{\eta!}
         \left|\partial^{\eta}
                 (\partial^{\gamma}u)\right|_{0,p,\bT}^p\right)} \notag \\
   & =  \frac{\left(\max_{i=1,2,3}\alpha_i\right)^{(k+1-m)p}
   \sum_{|\gamma| = m} \frac{m!}{\gamma!}(\alpha)^{- \gamma p}
         \left|\partial^{\gamma}u \right|_{0,p,\bT}^p}
       {\sum_{|\gamma| = m} \frac{m!}{\gamma!}
      (\alpha)^{-\gamma p}
        \left|\partial^{\gamma}u \right|_{k+1-m,p,\bT}^p} \notag \\
   & \le  \frac{\left(\max_{i=1,2,3}\alpha_i\right)^{(k+1-m)p}
  \sum_{|\gamma|=m} 
      \frac{m!}{\gamma!}(\alpha)^{- \gamma p}
              |\partial^{\gamma}u|_{0,p,\bT}^p }
       {\sum_{|\gamma|=m} \frac{m!}{\gamma!}(\alpha)^{- \gamma p}
         \left(A_p^{\gamma,k}\right)^{-1}
        |\partial^{\gamma}u|_{0,p,\bT}^p} \notag \\
  & \le C_{k,m,p}^p \left(\max_{i=1,2,3}\alpha_i\right)^{(k+1-m)p}, 
   \label{koyume}
\end{align}
}
where $C_{k,m,,p} := \max_{|\gamma|=m} A_p^{\gamma,k}$.
Here, we use the equality
\[
  \frac{(k+1)!}{\delta!} =
   \sum_{\substack{\gamma+\eta = \delta \\ |\gamma|=m, |\eta|=k+1-m}}
   \frac{m!}{\gamma!} \frac{(k+1-m)!}{\eta!}.
\]
Hence, Theorem~\ref{squeezingtheorem} is proved for this case. 
The proof of the case $p = \infty$ may be done in a similar manner.
 $\square$

\vspace{1mm}
\noindent
\textbf{Exercise:} (1) Check the above proof in detail.  For example,
confirm that, if $k=m=1$, \eqref{koyume} can be written as
{\allowdisplaybreaks
\begin{align*}
 & \frac{|v|_{1,p,T_{\alpha_1\alpha_2\alpha_3}}^p}
    {|v|_{2,p,T_{\alpha_1\alpha_2\alpha_3}}^p}
    = \frac{\sum_{|\gamma| = 1} \frac{1}{\gamma!}
      (\alpha)^{-\gamma p}\left|\partial^{\gamma}u\right|_{0,p,\bT}^p}
      { \sum_{|\delta|= 2} \frac{2!}{\delta!}
     (\alpha)^{- \delta p} \left|\partial^{\delta} u  \right|_{0,p,\bT}^p}
       \\
  & \quad = \frac{ \frac{1}{\alpha_1^p} |\partial_x u|_0^p +
             \frac{1}{\alpha_2^p} |\partial_y u|_0^p +
             \frac{1}{\alpha_3^p} |\partial_z u|_0^p}
       { \frac{1}{\alpha_1^{2p}} |\partial_{xx} u|_0^p +
         \frac{1}{\alpha_2^{2p}} |\partial_{yy} u|_0^p +
         \frac{1}{\alpha_3^{2p}} |\partial_{zz} u|_0^p +
         \frac{2}{\alpha_1^{p}\alpha_2^p} |\partial_{xy} u|_0^p +
         \frac{2}{\alpha_2^{p}\alpha_3^p} |\partial_{yz} u|_0^p +
         \frac{2}{\alpha_3^{p}\alpha_1^p} |\partial_{zx} u|_0^p }  \\
  & \quad = \frac{ \frac{1}{\alpha_1^p} |\partial_x u|_0^p +
             \frac{1}{\alpha_2^p} |\partial_y u|_0^p +
             \frac{1}{\alpha_3^p} |\partial_z u|_0^p}
           {\frac{1}{\alpha_1^p} X +
              \frac{1}{\alpha_2^p} Y +
             \frac{1}{\alpha_3^p} Z}  \\
 & \qquad {\scriptstyle
   \Bigl(X := \frac{1}{\alpha_1^p} |\partial_{xx} u|_0^p +
              \frac{1}{\alpha_2^p} |\partial_{xy} u|_0^p + 
              \frac{1}{\alpha_3^p} |\partial_{xz} u|_0^p, \quad
        Y := \frac{1}{\alpha_1^p} |\partial_{xy} u|_0^p +
              \frac{1}{\alpha_2^p} |\partial_{yy} u|_0^p + 
              \frac{1}{\alpha_3^p} |\partial_{yz} u|_0^p, \quad
    } \\ 
 & \hspace{6.7cm}  {\scriptstyle
    Z := \frac{1}{\alpha_1^p} |\partial_{zx} u|_0^p +
              \frac{1}{\alpha_2^p} |\partial_{zy} u|_0^p + 
              \frac{1}{\alpha_3^p} |\partial_{zz} u|_0^p
      \Bigr)} \\
   & \le  \frac{\left(\max_{i=1,2,3}\alpha_i\right)^{p}
    \left(\frac{1}{\alpha_1^p} |\partial_x u|_0^p +
             \frac{1}{\alpha_2^p} |\partial_y u|_0^p +
             \frac{1}{\alpha_3^p} |\partial_z u|_0^p\right)}
       {\frac{1}{\alpha_1^p} |\partial_x u|_1^p +
             \frac{1}{\alpha_2^p} |\partial_y u|_1^p +
             \frac{1}{\alpha_3^p} |\partial_z u|_1^p } \\
  & \qquad {\scriptstyle
   \Bigl(X \ge M |\partial_x u|_1^p, \quad
         Y \ge M |\partial_y u|_1^p, \quad
         Z \ge M |\partial_z u|_1^p, \quad
         M := \left(\max_{i=1,2,3}\alpha_i\right)^{-p}
    \Bigr)} \\ 
   & \le  \frac{\left(\max_{i=1,2,3}\alpha_i\right)^{p}
      \frac{(A_p^{(1,0,0),1})^p}{\alpha_1^p} |\partial_x u|_1^p +
             \frac{(A_p^{(0,1,0),1})^p}{\alpha_2^p} |\partial_y u|_1^p +
             \frac{(A_p^{(0,0,1),1})^p}{\alpha_3^p} |\partial_z u|_1^p}
       {\frac{1}{\alpha_1^p} |\partial_x u|_1^p +
             \frac{1}{\alpha_2^p} |\partial_y u|_1^p +
             \frac{1}{\alpha_3^p} |\partial_z u|_1^p}   \\
  & \le C_{1,1,p}^p \left(\max_{i=1,2,3}\alpha_i\right)^{p},
  \qquad C_{1,1,p} := \max\{A_p^{(1,0,0),1}, A_p^{(0,1,0),1}, A_p^{(0,0,1),1}\}. 
\end{align*}

\noindent
(2) Prove Theorem~\ref{squeezingtheorem} for the case $p = \infty$.

\end{document}